\documentclass[aop,seceqn,dvips]{arximspdf}
\usepackage{graphics}
\doi{10.1214/009117907000000213}
\volume{36}
\issue{1}
\pubyear{2008}
\firstpage{331}
\lastpage{362}

\makeatletter
\setattribute{keyword}{AMS}{AMS 2000 subject classification.}

\newproclaim{rem}{Remark}
\newtheorem{lem}{Lemma}
\newtheorem{cor}{Corollary}
\makeatother

\begin{document}
\begin{frontmatter}

\title{Shape fluctuations are different in different~directions}
\runtitle{Shape fluctuations are different in different directions}

\begin{aug}
\author[A]{\fnms{Yu} \snm{Zhang}\ead[label=e1]{yzhang3@uccs.edu}\thanksref{t1}\corref{}}
\thankstext{t1}{Supported by NSF Grant DMS-04-05150.}
\runauthor{Y. Zhang}
\affiliation{University of Colorado}
\address[A]{Department of Mathematics\\
University of Colorado\\
Colorado Springs, Colorado 80933\\
USA\\
\printead{e1}} %adresu isvedimo komanda gale!
\end{aug}

% HISTORY:
\received{\smonth{1} \syear{2006}}
\revised{\smonth{1} \syear{2007}}

% ABSTRACT
%
\begin{abstract}
We consider the first passage percolation model on $\mathbf{Z}^2$.
In this model, we assign independently to each edge $e$ a passage time $t(e)$
with a common distribution $F$. Let $T(u,v)$ be the passage time from
$u$ to $v$.
In this paper, we show that, whenever $F(0) < p_c$,
$\sigma^2( T((0,0),\break (n,0))) \geq C \log n $ for all $n\geq1$.
Note that if $F$ satisfies an additional special condition,
$\inf\operatorname{supp}(F)=r >0$ and $F(r) > \vec{p}_c$,
it is known that there exists $M$ such that for all $n$,
$\sigma^2(T((0,0), (n,n)))\leq M $.
These results tell us that shape fluctuations not only depend on distribution
$F$, but also
on direction.
When showing this result, we find the following interesting geometrical
property.
With the special distribution above,
any long piece with $r$-edges in an optimal path from $(0,0)$ to $(n,0)$
has to be very circuitous.
\end{abstract}

% KEYWORDS
%

\end{frontmatter}

%s1 ###
\section{Introduction of the model and results}\label{s1}

The first passage percolation model was introduced in 1965 by
Hammersley and
Welsh.
In this model, we consider the $\mathbf{Z}^2$ lattice as a graph with edges
connecting each pair of vertices
$u=(u_1, u_2)$ and $v=(v_1, v_2)$
with $d(u,v)=1$, where $d(u,v)$ is the Euclidean distance between $u$
and $v$.
We assign independently to each edge a nonnegative \textit{passage time} $t(e)$
with a common distribution $F$.
More formally, we consider the following probability space. As the
sample space, we take $\Omega=\prod_{e\in\mathbf{Z}^2} [0,\infty)$, whose points
are called \textit{configurations}.
Let $P=\prod_{e\in\mathbf{Z}^2} \mu_e$ be the corresponding product
measure on $\Omega$, where
$\mu_e$ is the measure on $[0, \infty)$ with distribution $F$. The
expectation and variance with respect to $P$ are denoted by $E(\cdot)$
and $\sigma^2(\cdot)$.
For any two vertices $u$ and $v$,
a path $\gamma$ from $u$ to $v$ is an alternating sequence
$(v_0, e_1, v_1,\ldots,v_i, e_{i+1}, v_{i+1},\ldots,v_{n-1},e_n, v_n)$ of vertices
$v_i$ and
edges $e_i$ between $v_i$ and $v_{i+1}$ in $\mathbf{Z}^2$ with $v_0=u$
and $v_n=v$.
Given such a path $\gamma$, we define its passage time as
%
%e1.1 ###
\begin{equation}\label{e1.1}
T(\gamma)= \sum_{i=1}^{n} t(e_i).
\end{equation}
For any two sets $A$ and $B$, we define the passage time from $A$ to $B$
as
\[
T(A,B)=\inf\{ T(\gamma)\},
\]
where the infimum is over all possible finite paths from some vertex in
$A$ to
some vertex in $B$.
A path $\gamma$ from $A$ to $B$ with $T(\gamma)=T(A, B)$ is called the
\textit{optimal path} of
$T(A, B)$. The existence of such an optimal path has been proven (see
Kesten~\cite{kes86}).
We also want to point out that the optimal path may not be unique.
If we focus on a special configuration $\omega$, we may write
$T(A,B)(\omega)$ instead of $T(A, B)$.
When $A=\{u\}$ and $B=\{v\}$ are single vertex sets, $T(u,v)$ is the passage
time
from $u$ to $v$. We may extend the passage time over $\mathbf{R}^2$.
If $x$ and $y$ are in $\mathbf{R}^2$, we define $T(x, y)=T(x', y')$, where
$x'$ (resp., $y'$) is the nearest neighbor of $x$ (resp., $y$) in
$\mathbf{Z}^2$. Possible indetermination can be eliminated by choosing an
order on the vertices of $\mathbf{Z}^2$ and taking the smallest nearest
neighbor for this order.

With these definitions, we would like to introduce the basic
developments and
questions
in this field. Hammersley and Welsh \cite{hamWel65} first studied the
point-point and
the point-line passage times
defined as follows:
\begin{eqnarray*}
a_{m,n} &=& \inf\{T(\gamma)\dvtx \gamma \mbox{ is a path from }
(m, 0)\mbox{ to }(n,0)\},
\\
b_{m,n} &=& \inf\bigl\{T(\gamma)\dvtx \gamma \mbox{ is a path from $(m, 0)$ to
$\{x=n\}$}\bigr\}.
\end{eqnarray*}
It is well known (see Smythe and Wierman \cite{smyWie78})
that if $Et(e) < \infty$,
%
%e1.2 ###
\begin{equation}\label{e1.2}
\lim_{n\rightarrow\infty}{1\over n} a_{0,n}
= \lim_{n\rightarrow\infty}{1\over n} b_{0,n}= \mu
\qquad \mbox{a.s. and in } L_1,
\end{equation}
where the nonrandom constant $\mu=\mu(F)$ is called the \textit{time constant}.
Later, Kesten showed (see Theorem 6.1 in Kesten \cite{kes86}) that
%
%e1.3 ###
\begin{equation}\label{e1.3}
\mu=0 \qquad \mbox{iff } F(0)\geq p_c,
\end{equation}
where $p_c=1/2$ is the critical probability for Bernoulli (bond)
percolation on $\mathbf{Z}^2$.

Given a vector $x\in\mathbf{R}^2$, by the same arguments as in
(\ref{e1.2}) and (\ref{e1.3}), if \mbox{$Et(e) < \infty$}, then
%
%e1.4 ###
\begin{eqnarray}\label{e1.4}
\lim_{n\rightarrow\infty}{1\over n} T(\mathbf{0}, nx)
&=& \inf_{n} {1\over n} E T(\mathbf{0}, nx)
\nonumber\\[-8pt]
\\[-8pt]
&=& \lim_{n\rightarrow\infty}
{1\over n} E T(\mathbf{0}, nx)=\mu(x) \qquad
\mbox{a.s. and in } L_1,
\nonumber
\end{eqnarray}
and
\[
\mu(x)=0\qquad \mbox{iff }F(0) \geq p_c.
\]
For convenience, we assume that $t(e)$ is not a constant and satisfies the
following:
%
%e1.5 ###
\begin{equation}\label{e1.5}
\int e^{\lambda x}\, dF(x) < \infty\qquad \mbox{for some } \lambda>0.
\end{equation}

When $F(0) < p_c$, the map $x \rightarrow\mu(x)$ induces a norm on
$\mathbf{R}^2$. The unit radius ball for this norm is denoted by
$\mathbf{B}:=\mathbf{B}(F)$
and is called the \textit{asymptotic shape}. The boundary of $\mathbf{B}$ is
\[
\partial\mathbf{B}:= \{ x \in\mathbf{R}^2\dvtx \mu(x)=1\}.
\]
$\mathbf{B}$ is a compact convex
deterministic set and $\partial\mathbf{B}$ is a continuous convex closed
curve (Kesten \cite{kes86}). Define for all $t> 0$,
\[
B(t):= \{v\in\mathbf{R}^2,  T( \mathbf{0}, v) \leq t\}.
\]
The shape theorem (see Theorem 1.7 of Kesten \cite{kes86}) is the
well-known result stating that for any $\varepsilon>0$,
\[
t\mathbf{B}(1-\varepsilon) \subset{B(t) } \subset t\mathbf{B}(1+\varepsilon)
\qquad \mbox{eventually w.p.1.}
\]
In addition to $t\mathbf{B}$, we can consider the mean of $B(t)$ to be
\[
G(t)=\{v\in\mathbf{R}^2\dvtx ET(\mathbf{0}, v)\leq t\}.
\]
By (\ref{e1.4}), we also have
\[
t\mathbf{B}(1-\varepsilon) \subset{G(t) } \subset t\mathbf{B}(1+\varepsilon).
\]

The natural and most challenging aspect in this field (see Kesten
\cite{kes86} and Smythe and Wierman \cite{smyWie78})
is to question the ``speed'' and ``roughness'' of the interface $B(t)$
from the deterministic boundaries $t\mathbf{B}$ and $G(t)$. This problem
has also
received a great amount of attention from
statistical physicists because of its equivalence with one version of
the Eden
growth model. They believe that there
is a scaling relation for the shape fluctuations in growth models. For each
unit vector $x$, we may denote by
$h_t(x)$ the \textit{height} of the interface (see page 490 in Krug and Spohn
\cite{kruSpo92}).
The initial condition is $h_0(x)=0$. Being interested in fluctuation, we
consider the height fluctuation function
\[
\bar{h}_t(x)= h_t(x)- Eh_t(x).
\]

Statistical physicists believe that $\bar{h}_t(x)$ should satisfy
(see (3.1) in Krug and Spohn \cite{kruSpo92}) the following scaling property
\[
\bar{h}_t(x)= b^{\zeta}\bar{h}_{b^{z}t}(bx)
\]
with the scaling exponents $\zeta$ and $z$ for an arbitrary rescaling
factor $b$.
With this scaling equation, we should have (see (7.9) in Krug and Spohn
\cite{kruSpo92}), for all vectors~$x$,
%
%e1.6 ###
\begin{equation}\label{e1.6}
\bar{h} _t(x)\approx t^{\zeta/ z} \qquad \mbox{pointwisely} \quad \mbox{or}\quad
\sigma(h_t(x))\approx t^{\zeta/ z}.
\end{equation}
In particular, it is believed that $\zeta=1/2$ and $z=2/3$ when $d=2$.

Mathematicians have also made significant efforts in this direction.
When $F(0) > p_c$, it is known (see Zhang \cite{zha05}) that
%
%e1.7 ###
\begin{equation}\label{e1.7}
\sigma^2 (a_{0,n})< \infty.
\end{equation}

When $F(0)=p_c$ and $t(e)$ only takes two values zero and one, it is
also known
(see Kesten and Zhang \cite{kesZha97}) that
%
%e1.8 ###
\begin{equation}\label{e1.8}
\sigma^2(a_{0,n})=O(\log n).
\end{equation}
In fact, Kesten and Zhang \cite{kesZha97} showed a CLT for the process
$a_{0,n}$, a much stronger
result than (\ref{e1.8}). For a more general distribution $F$ with $F(0)=p_c$,
$\sigma^2(a_{0,n})$
can be either convergent or divergent, depending on the behavior of the
derivative
of $F(x)$ at $x=0$ (see Zhang \cite{zha99}).

Now we focus on the most interesting situation: when $F(0) < p_c$.
It is widely conjectured (see (\ref{e1.6}) above and Kesten
\cite{kes93}) that if $F(0) < p_c$, then
\[
\sigma^2 (a_{0,n})\approx n^{2/3}.
\]
The mathematical estimates for the upper bound of $\sigma^2(a_{0,n})$
are quite promising. Kesten \cite{kes93}
showed that if $F(0)< p_c$, there is a constant $C_1$ such that
%
%e1.9 ###
\begin{equation}\label{e1.9}
\sigma^2( a_{0,n})\leq C_1 n.
\end{equation}
In this paper, $C$ and $C_i$ are always positive constants that may
depend on
$F$, but not on $t$, $m$, or $n$.
Their values
are not significant and change from
appearance to appearance.
Benjamini, Kalai and Schramm \cite{benKalSch03}
also showed that when $t(e)$ only takes two
values $0<a< b$ with a half probability for each one,
\[
\sigma^2( a_{0,n})\leq C_1 n/\log n,
\]
where $\log$ denotes the natural logarithm.

\begin{figure}

\includegraphics{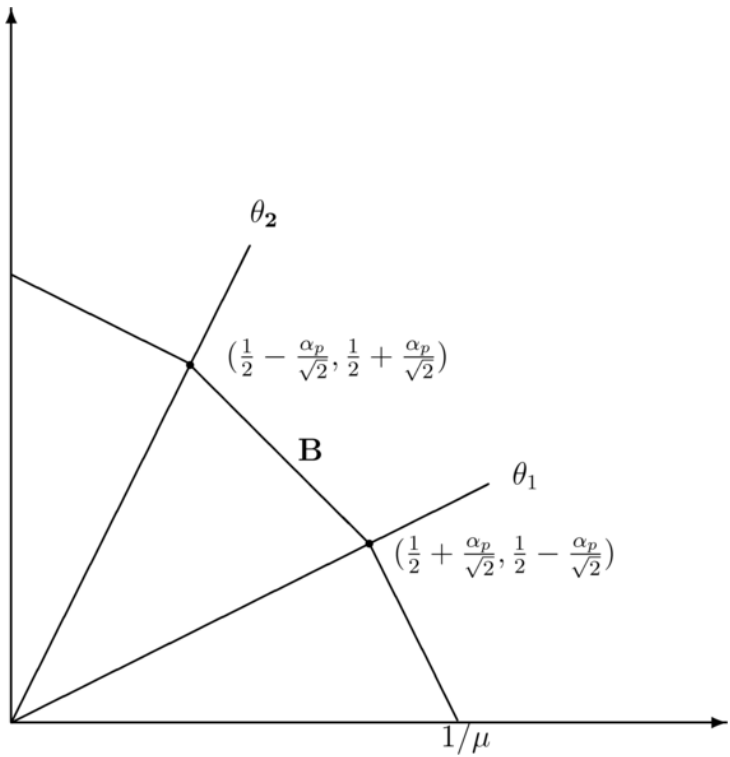}

\caption{The graph shows that shape $\mathbf{B}$ contains a flat
segment when $F$ satisfies (\protect\ref{e1.10}).\label{f1}}
\end{figure}

On the other hand, the lower bound of the variance for $\sigma^2(a_{0,n})$
seems to be much more difficult to estimate.
For a high-dimensional lattice, there are some discussions for a lower
bound of the fluctuations from
$B(t)$ to $t\mathbf{B}$ (see Zhang \cite{zha06}).
In this paper, we would like to focus
on the square lattice.
To understand the complexity of the lower bound, we have to deal with the
following special distributions
investigated by Durrett and Liggett \cite{durLig81}. They defined
\[
r=\inf\operatorname{supp}(F)=\inf\bigl\{x \dvtx F(x)
=P\bigl(t(e) \leq x\bigr)>0\bigr\}
\]
with $r >0$.
Clearly, if $r >0$, $F(0)=0 < p_c$, so shape $\mathbf{B}$ is compact.
Without loss of generality, we can suppose that $r=1$ if we replace
$F(x)$ by $F(rx)$. In the following, we always assume that
%
%e1.10 ###
\begin{equation}\label{e1.10}
\inf\operatorname{supp} (F)=1 \quad \mbox{and}\quad  F(1)
= P\bigl(t(e) =1\bigr) \geq \vec{p}_c,
\end{equation}
where $\vec{p}_c$ is the critical value for the oriented percolation model.
Durrett and Liggett \cite{durLig81}
found that shape $\mathbf{B}$ contains a flat segment
on the diagonal direction.
Later, Marchand \cite{mar02} presented the precise locations of
the flat segment in the shape when distribution $F$ satisfies (\ref{e1.10}).
More precisely, two polar coordinates in the first quadrant are denoted by
$(\sqrt{1/2+\alpha^2_p}, \theta_i)$ for $i=1,2$ (see Figure~\ref{f1}), where
\[
\theta_i= \operatorname{arc}\tan
\biggl( { 1/2\mp\alpha_p/\sqrt{2} \over 1/2\pm\alpha_p/\sqrt{2} } \biggr),
\]
and $\alpha_p\geq0$ is a constant defined in (2.4) below.
Note that $\theta_1 < \theta_2$ if $F(1) > \vec{p}_c$ and
$\theta_1=\theta_2$ if $F(1)=\vec{p_c}$.
Marchand (see Theorem 1.3 of Marchand \cite{mar02})
showed that under (\ref{e1.10}),
\begin{eqnarray*}
&& \mathbf{B}\cap\{(x,y)\in\mathbf{R}^2, |x|+|y|=1\}
\\
&&\qquad = \mbox{the segment from }
\bigl(\sqrt{1/2+\alpha^2_p}, \theta_1\bigr)
\mbox{ to } \bigl(\sqrt{1/2+\alpha^2_p},\theta_2\bigr),
\end{eqnarray*}
where the segment will shrink as a point $(1/\sqrt{2}, \pi/4)$ when
$F(1)=\vec{p}_c$.
This segment is called the \textit{flat edge} of shape $\mathbf{B}$.
The cone between $\theta_1$ and $\theta_2$ is called
the \textit{oriented percolation cone}.
To understand why this is called the oriented percolation cone, we
introduce the
following oriented paths. Let us define the northeast-
and the southeast-oriented paths.
A path (not necessary to be a $1$-path) is said to be a northeast
path if each vertex $u$ of the path has only one exiting edge,
either from $u$ to $u+(1,0)$ or to $u+(0,1)$.
Similarly, a path is said to be a southeast path if each vertex $u$ of the
path has only one existing edge,
either from $u$ to $u+(1,0)$ or to $u+(0,-1)$.

%-------------------------
%f1 ###
%%-------------------------

For any vector $(r, \theta)$ with $\theta_1 \leq\theta\leq\theta_2$, under
(\ref{e1.10}), with a positive probability,
there is a northeast path $\gamma$ from $(0,0)$ to $(nr, \theta)$ with
only 1-edges
(see (3.2) in Yukich and Zhang \cite{yukZha06}). Thus, we call the cone between
$\theta_1$ and $\theta_2$
the oriented percolation cone.

With this observation, for
$\theta_1 < \theta< \theta_2$ (see special case 1 in Newman and Piza~\cite{newPiz95}),
there exists a constant $C=C(F, \theta)$ such that
%
%e1.11 ###
\begin{equation}\label{e1.11}
\sigma^2(T((0,0), (n, \theta))) \leq C \qquad \mbox{for all }n,
\end{equation}
where both $(0,0)$ and $(n, \theta)$ are polar coordinates.
On the other hand, it has been proven (see Newman and Piza (1995)) that if
$F(0) <p_c$ and $\inf\operatorname{supp}(F)=0$, or
$\inf\operatorname{supp}(F)=1$ and $F(1)< \vec{p}_c$, then
%
%e1.12 ###
\begin{equation}\label{e1.12}
\sigma^2(T((0,0), (n, \theta)))\geq C \log n
\qquad \mbox{for all }\theta.
\end{equation}
Even though (\ref{e1.12}) is far from the correct order $n^{2/3}$, it at least
tells us that
$\sigma^2(a_{0,n})$ diverges as $n\rightarrow\infty$. As we mentioned earlier,
both the convergence and divergence in
(\ref{e1.11}) and (\ref{e1.12}) indicate the complexity
of an estimate for the lower bound of the variance.

From (\ref{e1.11}) and (\ref{e1.12}), we may ask the behaviors of the variance for
the passage time on a nonoriented percolation cone or simply ask
whether with (\ref{e1.10}), for the most popular first passage time $a_{0,n}$,
%
%e1.13 ###
\begin{equation}\label{e1.13}
\sigma^2(a_{0,n}) \qquad \mbox{diverges as }n\rightarrow\infty.
\end{equation}
Indeed, if there were a proof for (\ref{e1.13}), the proof would be tricky
because one has to show that
two different behaviors exist in the oriented percolation and the nonoriented
percolation cones.
As Newman and Piza (1995) described, ``either the new techniques, or additional
hypotheses seem to need to investigate conjecture~(\ref{e1.13})
when (\ref{e1.10}) holds.'' This is the same flavor as the
extension of the strict inequality on the time constant (see van den
Berg and Kesten \cite{benKalSch03})
to the nonoriented percolation cone (see Marchand \cite{mar02}).
In this paper, one of the main works is to investigate the different
behaviors in the oriented and nonoriented percolation cones. We
discovered that,
unlike the oriented
percolation cone, any long piece with $1$-edges in an optimal path from
$(0,0)$ to $(n,0)$
contains proportional circuitous pieces. With this geometric property,
we will show (\ref{e1.13}).

Before we mention our result, we would like to introduce Newman and Piza's
martingale method that was used to show
(\ref{e1.12}). To describe their method simply, we assume that $t(e)$ can only take
two values $a$ and $b$ with $0\leq a < b$.
The key to these martingale arguments is to show that there are proportionally
many $b$-edges
in an optimal path. If we change $b$-edges from $b$ to $a$, then
the passage time will shrink at least $b-a$ from the original passage time.
This tells us why the variance of passage time should be
large if the number of $b$-edges is large. The remaining task is to
estimate the number of $b$-edges.
However, if we assume (\ref{e1.10}) holds, it is possible that all edges in any
optimal path have value $a$.
Therefore, the Newman and Piza method will not be applied.

%-----------------------
%f2 ###
\begin{figure}[b]

\includegraphics{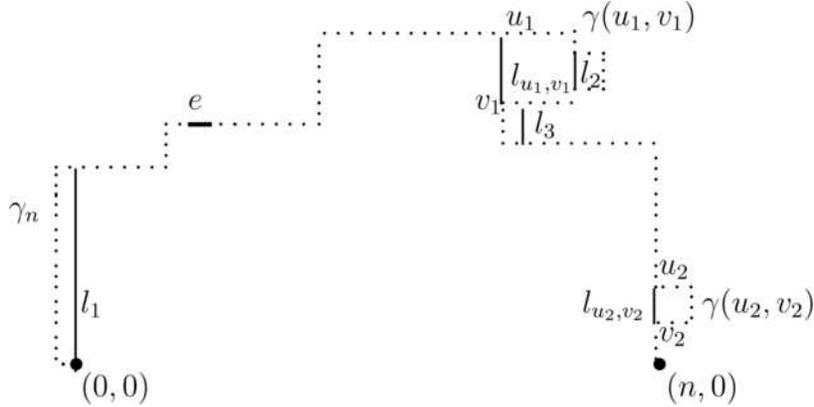}

\caption{The dotted line in the graph is the optimal path $\gamma_n$.
$e$ is a $1^+$-edge.
$l_1$, $l_2$, $l_3$, $l_{u_1,v_2}$ and $l_{u_2, v_2}$ are broken
bridges of $\gamma_n$.
In fact, there are many vertical $M$-broken bridges
parallel to $l_{u_1, v_1}$ and $l_3$, but we will not list all of them
in this graph.
$l_{u_2, v_2}$ is the only broken bridge for the remaining path from
$v_1$ to $(n,0)$.
$l_1$ is not an $M$-bridge because its length is more than $2M$.
After ordering all $M$-broken bridges of $\gamma_n$, $ l_2$ and $l_3$
are no longer $M$-broken bridges
for $\gamma_n$, so $l_{u_1, v_1}$ and $l_{u_2, v_2}$ are $M$-broken
bridges for $\gamma_n$.\label{f2}}
\end{figure}
%%-----------------------

To develop new techniques for case (\ref{e1.10}),
we need to investigate the geometric properties of an optimal path.
For an optimal path from $(0, 0)$ to $(n,0)$, we may guess that an
optimal path should not be northeast or southeast since it is not
in the oriented percolation cone.
Let us give a more precise definition of what
``nonnortheast'' or ``nonsoutheast'' means.
Let $\gamma_n$ be an optimal path from $(0,0)$ to $(n,0)$.
Note that the existence of such a $\gamma_n$ has been mentioned before.
With this existence, there might be many such optimal paths for $a_{0,n}$.
We now select a unique optimal path. For example, we may start at the
origin and select vertices among all optimal paths
in each step such that
the selected vertex is closer
to the $X$-axis. We still write this unique optimal path as
$\gamma_n$, without loss of generality. Later, we always use $\gamma_n$
as the
unique selected optimal path.
For vertices $u,v$, we say $l_{u,v}$ is an \textit{$M$-bridge}
if $l_{u,v}$ is a horizontal or vertical segment, including both
vertices and edges, from $u$ to $v$ whose number of vertices is
less than $2M$.
Furthermore, we say $l_{u,v}$ is an $M$-\textit{broken bridge} of
$\gamma_n$ if $l_{u,v}$ is an $M$-bridge and
%
%e1.14 ###
\begin{equation}\label{e1.14}
l_{u,v}\cap\gamma_n=\{u,v\}.
\end{equation}
In other words, the bridge from $u$ to $v$ is broken and $\gamma_n$ has
to go
around from $u$ to $v$ to
avoid using any vertex in $l_{u,v}$ except $u$ and $v$.

Now we choose a list of peculiar broken bridges of $\gamma_n$ (see
Figure~\ref{f2}), as follows.
We first list all the possible $M$-broken bridges of $\gamma_n$.
We then go along $\gamma_n$ from the origin to meet $u_1$,
the first vertex of $\gamma_n$,
such that there exists an $M$-broken bridge for $\gamma_n$ at $u_1$.
Note that there may be up to three such $M$-broken bridges at $u_1$, since
$u_1$ may be the origin.
If there are two or three $M$-broken bridges at $u_1$, for example,
$l_{u_1,v_1'}$, $l_{u_1,v_1''}$,
and $l_{u_1,v_1'''}$, we use the following way to select one.
Going along $\gamma_n$ from the origin, we first meet either
$v_1'$, $v_1''$ or $v_1'''$. We then select the first vertex that we encounter
and denote it by $v_1$.
We can go along $\gamma_n$ from $u_1$ to $v_1\neq u_1$, denoted by
$\gamma(u_1,v_1)$.
Note that $\gamma(u_1, v_1)\cup l_{u_1,v_1}$ is a loop (see Figure~\ref{f2}).
After selecting $l_{u_1, v_1}$, we list all the possible $M$-broken
bridges of
the path from $v_1$ to $(n,0)$
along the remaining part of $\gamma_n$.
As we go from $v_1$ along the remaining part of $\gamma_n$, we meet $u_2$,
the first vertex in the remaining part, such that there exists an $M$-broken
bridge at $u_2$ for
the remaining part. Note that $u_2$ may equal $v_1$.
Similarly, if there are two $M$-broken bridges at $u_2$, denoted by
$l_{u_2,v_2'}$
and $l_{u_2,v_2''}$, we select $v_2$, from $v_2'$ and $v_2''$,
as the first vertex encountered on the remaining part of
$\gamma_n$ from $v_1$.
We now go along the remaining part of $\gamma_n$ from $u_2$ to $v_2$,
denoted by
$\gamma(u_2, v_2)$.
Thus, $\gamma(u_2, v_2)\cup l_{u_2,v_2}$ is the second loop.
Since $\gamma_n$ is finite, we continue this process until the last $M$-broken
bridge, $l_{u_\tau, v_\tau}$.
The corresponding piece of $\gamma_n$ from $u_\tau$ to $v_\tau$ is
$\gamma(u_\tau, v_\tau)$, and the
loop is $\gamma(u_\tau, v_\tau)\cup l_{u_\tau, v_\tau}$.
In the following discussion, for $\gamma_n$,
we always consider these $M$-broken bridges $\{l_{u_i, v_i}\}$
($i=1,\ldots,\tau$) for $\gamma_n$ by
this arrangement and by ignoring the other listed $M$-broken bridges.

Furthermore, by the definition (see Figure~\ref{f2}),
%
%e1.15 ###
\begin{eqnarray}\label{e1.15}
&& \mbox{the subpath of $\gamma_n$ from $v_i$ to $u_{i+1}$ has none}
\nonumber\\[-8pt]
\\[-8pt]
&&\qquad \mbox{of its own $M$-broken bridges.}
\nonumber
\end{eqnarray}
Note that $\gamma_n$ is self-avoiding, so
%
%e1.16 ###
\begin{equation}\label{e1.16}
\gamma(u_i, v_i)\quad \mbox{and}\quad  \gamma(u_j, v_j)
\mbox{ have no common edge for }i\neq j .
\end{equation}
On the other hand,
the interior of the loop $\gamma(u_i, v_i)\cup l_{u_i, v_i}$ cannot
contain a vertex of $\gamma_n$
by (\ref{e1.16}) and the definition of the $M$-broken bridge. With this observation,
%
%e1.17 ###
\begin{equation}\label{e1.17}
l_{u_i, v_i}\mbox{ and } l_{u_j,v_j}\mbox{ have no common edge for }
i\neq j.
\end{equation}
Since $l_{u_i, v_i}$ is shorter than $\gamma(u_i, v_i)$ by at least two edges,
and each edge costs at least
time one,
%
%e1.18 ###
\begin{equation}\label{e1.18}
\exists  e\in l_{u_i, v_i}\mbox{ such that } t(e) >1.
\end{equation}

Clearly, for a northeast or southeast path, there is no broken bridge.
From this point of view,
we may guess that there are many broken bridges for the optimal path
$\gamma_n$ from the origin to $(n,0)$.
We shall show the following theorem to describe this fact.

An edge $e$ is called a 1-edge if $t(e)=1$. A path is called a 1-path
if all
of its edges are 1-edges.
Note that we assume that $t(e)$ is not a constant, so
%
%e1.19 ###
\begin{equation}\label{e1.19}
0< P\bigl(1<t(e)\bigr).
\end{equation}
We say edge $e\in\gamma_n$ is a $1^+$-edge if $t(e)>1$. We collect all vertices
in $\gamma_n$
that are adjacent to $1^+$-edges on $\gamma_n$ and denote them by
$D(\gamma_n)$.
If $\gamma_n$ is not northeast or southeast, there may exist $M$-broken bridges
$\{l_{u_i,v_i}\}_{1\leq i\leq\tau}$ of $\gamma_n$. Note that
\mbox{$u_i, v_i\in \gamma_n$}, so
we collect all vertices $u_i$ and $v_i$ in $\gamma_n$ for $0\leq i\leq \tau$ and
denote them by $S_M(\gamma_n)$.
With these definitions, we will have the following theorem.
\begin{thm}\label{thm1}%thm1
Let $F$ be a distribution such that
$\inf\operatorname{supp}(F)=1$, $\vec{p}_c \leq F(1)$
and satisfying the tail assumption
in (\ref{e1.5}). Then there exist
constants $\delta=\delta(F, M)>0$ and $C_i=C_i(F, M,\delta)$ $(i=1,2)$
such that
for all $m,n\geq1$ with $n/2\geq m \geq n^{2/3}$,
\[
P\bigl(| B(m)\cap[S_M(\gamma_n)\cup D(\gamma_n)]| \leq\delta m\bigr)
\leq C_1\exp(-C_2 n^{1/14}),
\]
where $B(m)=[-m,m]^2$ and $|A|$ represents the number of vertices in
set $A$.
\end{thm}
\begin{rem}\label{rem1}%rem1
When $m=O(n)$, we can generalize Theorem \ref{thm1} for any vector $x=(1,\theta)$ with
$0< \theta< \theta_1$. More precisely, for a polar coordinate
$x=(1, \theta)$ with $0< \theta< \theta_1$, let $\gamma_n(\theta)$
be an optimal path from the origin to $(n, \theta)$. Similarly, we
choose a list of peculiar
$M$-broken bridges $\{ l_{u_i, v_i}(\theta)\}$
for $\gamma_n(\theta)$ as we did for $\gamma_n$.
We denote by $D(\gamma_n(\theta))$ all vertices in $\gamma_n(\theta)$
that are adjacent to $1^+$-edges on $\gamma_n(\theta)$.
We also denote by $S_M(\gamma_n(\theta))$ all vertices $\{u_i, v_i\}$
for the $M$-broken bridges $\{l_{u_i,v_i(\theta)}\}$ of $\gamma_n(\theta)$.
If $m=O(n)$, we can show, under the assumptions of Theorem~\ref{thm1},
%
%e1.20 ###
\begin{equation}\label{e1.20}
P\bigl(| B(m)\cap\bigl(S_M(\gamma_n(\theta))\cup D(\gamma_n(\theta))\bigr)| \leq\delta
n\bigr)\leq C_1\exp(-C_2 n).
\end{equation}
However, because of the lack of symmetry, we cannot show (\ref{e1.20})
for all $\theta< \theta_1$ when $m=o(n)$.
\end{rem}
\begin{rem}\label{rem2}
We may consider the same problem as Theorem \ref{thm1} when \mbox{$d\geq3$}.
It is possible to show a similar result as Theorem \ref{thm1}
when $d\geq3$ and (\ref{e1.10}) holds. However, we do not know
whether a similar result of Theorem 1 holds when
$\vec{p}_c(d) \leq F(1) < \vec{p}_c$,
where $\vec{p}_c(d)$ is a critical probability for the $d$-dimensional oriented
percolation. The main reason is that we need to use Lemma \ref{lem3}, proven by Marchand
\cite{mar02}, in our Section~\ref{s2}
to show Theorem 1, but Lemma \ref{lem3} has not been proven for all $d \geq3$.
\end{rem}
\begin{rem}\label{rem3}
The term $n^{2/3}$ in Theorem \ref{thm1} can be improved to $Cn^{1/2}
\log n$ for large constant $C$.
\end{rem}

With Theorem \ref{thm1}, we can see that an optimal
path contains proportionally many $1^+$-edges or
proportionally many vertices adjacent to $M$-broken bridges.
If we change $1^+$ edge in $\gamma_n$ to $1$-edge, or recover the
bridge by
changing the time of the $1^+$-edges
from $1^+$ to 1, we have saved a positive
passage time for $\gamma_n$.
Therefore, we can also use Newman and Piza's
\cite{newPiz95} martingale method, but with
a large square construction, to show the following theorem.
\begin{thm}\label{thm2}
Let $F$ be a distribution such that
$\inf\operatorname{supp}(F)=1$, $\vec{p}_c \leq F(1)$ and satisfying the tail assumption
in (\ref{e1.5}). Then there exists $C=C(F)$ such that for all $n$
\[
\sigma^2(a_{0,n})\geq C \log n .
\]
\end{thm}
\begin{rem}\label{rem4}
Together with Newman and Piza's result \cite{newPiz95} in
(\ref{e1.12}), we have
\[
\sigma^2(a_{0,n})\geq C \log n
\]
for all $n$ whenever $F(0) < p_c$.
Together with (\ref{e1.7}), (\ref{e1.8}) and (1.21),
the whole picture of convergence or
divergence for $\sigma^2(a_{0,n})$ is complete.
\end{rem}
\begin{rem}\label{rem5}
Under the same assumptions as in Theorem \ref{thm2}, the same
proof can be carried out to show that
\[
\sigma^2(b_{0,n}) \geq C \log n.
\]
\end{rem}
\begin{rem}\label{rem6}
We are unable to show Theorem \ref{thm2} for the passage time $T((0,0),
(n, \theta))$ for all $0< \theta< \theta_1$
even though we believe it is true. In fact, if one can show (\ref{e1.20}) in
Remark \ref{rem6} for all $n^{2/3}\leq m\leq n/2$,
then the same proof of Theorem \ref{thm2} can be carried out to show
\[
\sigma(T((0,0), (n, \theta))) \geq C\log n
\]
for all $0< \theta< \theta_1$.
\end{rem}
\begin{rem}\label{rem7}
As we mentioned in (\ref{e1.11}), there exists $C=C(F, \theta)$ for
$\theta_1< \theta< \theta_2$,
\[
\sigma^2 (T((0,0), (n, \theta))) \leq C \qquad \mbox{for all }n.
\]
This result can be generalized for $\theta=\theta_1$ and $\theta=\theta_2$
without too many difficulties.
\end{rem}

%s2 ###
\section[Preliminaries for the proof of Theorem 1]{Preliminaries for the proof of Theorem \protect\ref{thm1}}\label{s2}

Before presenting the proofs of the theorems we would like to introduce
a few lemmas.
\begin{lem}\label{lem1}
If $\gamma$ is a path with $|\gamma|\leq2M$ and without
$M$-broken bridge,
then $\gamma$ is either northeast or southeast.
\end{lem}
\begin{pf}
Denote by $u=(u_1, u_2)$ and $v=(v_1, v_2)$ the
extremities of $\gamma$.
By symmetry, we can assume that $v_1 \geq u_1$ and $v_2\geq u_2$. Let
us show that in this case $\gamma$ is northeast.

If $\gamma$ is not northeast, there exist, along $\gamma$, two successive
vertices $x=(x_1, x_2)$ and $y=(y_1, y_2)$ such that:
\begin{itemize}
\item either $y_1=x_1-1$ and $y_2=x_2$,
\item or $y_1=x_1$ and $y_2=x_2-1$.
\end{itemize}
Let us consider the first case. By the continuity of the first
coordinate along $\gamma$,
since $\gamma$ is simple, there necessarily exists
$k\in\{-2M+1,-2M+2,\ldots, -1,1,\ldots, 2M-2, 2M-1\}$
such that $y'=(y_1, y_2+k)$ and $x'=(x_1, x_2+k)$ are successive
vertices along $\gamma$. But now two cases occur:
\begin{itemize}
\item If the couple $(x, y)$ appears in $\gamma$ before the couple $(y',x')$, then
the vertical segment between $x$ and $x'$ contains an $M$-broken bridge for
$\gamma$.
\item If the pair $(y', x')$ appears in $\gamma$ before the pair
$(x, y)$, then the vertical segment between $y$, and $y'$ contains an $M$-broken
bridge for $\gamma$.
The second case can be proven by a similar argument.\quad \qed
\end{itemize}\noqed
\end{pf}

If we rotate our lattice counterclockwise by $45^\circ$ and extend each
edge by a factor of $\sqrt{2}$, the new graph is
denoted by $\mathcal{L}$ with
oriented edges from $(m,n)$ to $(m+1, n+1)$ and to $(m-1, n+1)$.
Each edge is independently open or closed with probability $p=F(1)$ or $1-p$.
For two vertices $u$ and $v$ in $\mathcal{L}$, we say $u\rightarrow v$ if
there is a sequence $v_0=u, v_1,\ldots, v_m=v$ of points of $\mathcal{L}$
with the vertices $v_i=(x_i,y_i)$ and $v_{i+1}=(x_{i+1}, y_{i+1})$ for
$0\leq i\leq m-1$ such that $y_{i+1} =y_i+1$ and $v_i$
and $v_{i+1}$ are connected by an open edge.
For $A\subset(-\infty, \infty)$, we denote a random subset by
\[
\xi_n^A=\{x\dvtx \exists x'\in A \mbox{ such that } (x',0)
\rightarrow (x,n) \} \qquad \mbox{for } n>0.
\]
The right edge for this set is defined by
\[
r_n=\sup\xi_n^{(-\infty, 0]}  \qquad  (\sup\varnothing=-\infty).
\]
We know (see Section 3 (7) in Durrett \cite{dur84}) that
there exists a nonrandom constant $\alpha_p$ such that
%
%e2.1 ###
\begin{equation}\label{e2.1}
\lim_{n\rightarrow\infty} {r_n\over n}=\lim_{n}\frac{E r_n}{ n}
=\alpha_p \qquad \mbox{a.s. and in }L_1,
\end{equation}
where $\alpha_p >0$ if $p> \vec{p}_c$ and $\alpha_p=0$ if $p=\vec{p}_c$.
Now we need to investigate the large deviation for the upper tail of $r_n$.
The exponential bound, when $p=F(1) > \vec{p}_c$, has been obtained by Durrett
\cite{dur84} in his Section 11.
However, his proof will not apply for $p =\vec{p}_c$.
We present a new proof, also independently interesting, to cover the case
$p=\vec{p}_c$.
\begin{lem}\label{lem2}
If $p=F(1)\geq\vec{p}_c$, for every $\eta>0$, there exist
$C_i=C_i(p,\eta)$ for $i=1,2$ such that
%
%e2.2 ###
\begin{equation}\label{e2.2}
P\bigl(r_n \geq n(\alpha_p +\eta)\bigr)\leq C_1 \exp(-C_2 n).
\end{equation}
\end{lem}
\begin{pf}
We observe that $r_n$ can be embedded in a two-parameter process (see Section 3
in Durrett \cite{dur84}). For $0\leq m< n$, let
\[
{r}_{m,n}=\sup\{x-r_m\dvtx (x, n)\in\mathcal{L}
\mbox{ and }
\exists y\leq r_m \mbox{ such that }(y, m)\rightarrow(x,n)\}.
\]
In particular, we denote by
\[
{r}_{m,n}(j)=\sup\{x-j\dvtx (x, n)\in\mathcal{L}
\mbox{ and }
\exists y\leq j \mbox{ such that }(y, m)\rightarrow(x,n)\}.
\]
It follows from Section 3 (3) and (4) in Durrett \cite{dur84} that
%
%e2.3 ###
\begin{equation}\label{e2.3}
r_{m, n}\stackrel{d}{=} r_{n-m}
\quad \mbox{and}\quad
r_n\leq r_m + r_{m,n} \qquad \mbox{for } 0\leq m < n.
\end{equation}
By (\ref{e2.1}), we take $M$ such that
%
%e2.4 ###
\begin{equation}\label{e2.4}
Er_M \leq M(\alpha_p +\eta/2).
\end{equation}
Without loss of generality, we may assume that $n/M=l$ is an integer.
By (\ref{e2.4}), we have
%
%e2.5 ###
\begin{equation}\label{e2.5}
P\bigl(r_n \geq n(\alpha_p+\eta) \bigr)\leq P ( r_n-l Er_M\geq n\eta/2 ) .
\end{equation}
By (\ref{e2.5}) and Markov's inequality, for any $t>0$,
%
%e2.6 ###
\begin{equation}\label{e2.6}
P\bigl(r_n \geq n(\alpha_p+\eta) \bigr)
\leq\exp(-tn\eta/2) E\exp [t (r_n - lEr_M ) ].
\end{equation}
By (\ref{e2.3}),
\[
E\exp [t ( r_n -lEr_M ) ]\leq E\exp [t  ( r_{n-M}
+r_{n-M,n} -l Er_M )  ].
\]

By our definition, $r_{n-M}$ only depends on the open and closed edges
in the region between $\{y=0\}$ and $\{y=n-M\}$. On the other hand, on
$\{r_{n-M}=j\}$
for some~$j$, $r_{n-M,n}=r_{n-M, n}(j)$ only depends on the open and closed
edges in the
region between $\{y=n-M\}$ and $\{y=n\}$. In addition, for any $j$,
\[
r_{n-M,n}(j)\stackrel{d}{=}r_M.
\]
With these observations,
\begin{eqnarray*}
&& E\exp [t ( r_n -lEr_M ) ]
\\
&&\qquad \leq E\exp [ t  ( r_{n-M} +r_{n-M,n} -l Er_M  ) ]
\\
&&\qquad = \sum_{j} E\exp
\bigl[t  \bigl( r_{n-M} -({l-1}) Er_M +r_{n-M, n}- Er_M \bigr) \bigr]I(r_{n-M}=j)
\\
&&\qquad = \sum_{j} E\exp
\bigl[t  \bigl(r_{n-M} -({l-1}) Er_M \bigr) \bigr]
I(r_{n- M}=j)\exp\bigl[t\bigl(r_{n-M, n}(j)-Er_M\bigr)\bigr]
\\
&&\qquad = \bigl\{E\exp \bigl[t  \bigl(r_{n-M} -({l-1}) Er_M \bigr) \bigr]
\bigr\}\{E\exp[t(r_M-Er_M)]\},
\end{eqnarray*}
where $I(A)$ is the indicator for event $A$.
We iterate this way $l$ times to have
%
%e2.7 ###
\begin{equation}\label{e2.7}
E\exp [ t  ( r_n -l Er_M ) ]\leq[E\exp t (r_M- Er_M)]^{l}.
\end{equation}
Note that $(r_M-Er_M)<\infty$ almost surely, so we use Taylor's
expansion for $\exp t (r_M-Er_M)$ to have
%
%e2.8 ###
\begin{equation}\label{e2.8}
E\exp t (r_M-Er_M)=E\sum_{i=0}^{\infty} {[t(r_M-Er_M)]^i\over i!}.
\end{equation}
If we can show that
%
%e2.9 ###
\begin{equation}\label{e2.9}
\qquad
\sum_{i=2}^\infty t^{i-2}{E|r_M-Er_M|^i\over i!}
\qquad \mbox{converges uniformly for $t\in[0,1/(4M)]$,}
\end{equation}
then by (\ref{e2.8}), for all $t\in[0,1/(4M)]$,
%
%e2.10 ###
\begin{eqnarray}\label{e2.10}
E\exp t (r_M-Er_M)
&\leq&  1+ t^2\sum_{i=2}^{\infty} t^{i-2}{E|(r_M- Er_M)|^i\over i!}
\nonumber\\[-8pt]
\\[-8pt]
&\leq & \bigl(1+O(t^2)\bigr)\leq\exp(C t^2)
\nonumber
\end{eqnarray}
for some constant $C=C(M, F)$.
By (\ref{e2.6}), (\ref{e2.7}) and (\ref{e2.10}),
if we take $t\leq1/(4M)$ small enough, then there
exist $C_i=C_i(F, \eta)$ for $i=1,2$, such that
%
%e2.11 ###
\begin{equation}\label{e2.11}
P\bigl(r_n \geq n(\alpha_p+\eta)\bigr) \leq C_1 \exp(-C_2 n),
\end{equation}
so Lemma 2 follows. Now it remains to show (\ref{e2.9}).
Note that $r_1$ is of geometric type with $p\geq\vec{p}_c >0$, so by (\ref{e2.3}),
there exists $C_3(M, F)$ such that
%
%e2.12 ###
\begin{equation}\label{e2.12}
E|r_M|^i \leq E(M|r_1|)^i\leq C_3M^ii! .
\end{equation}
By (\ref{e2.12}), there exists a constant $C_4=C_4(M, F)$
such that
%
%e2.13 ###
\begin{eqnarray}\label{e2.13}
\qquad \hspace*{12pt}
E|r_M-Er_M|^i &\leq& E(|r_M| +E|r_M|)^i
\nonumber\\
&=& E[(|r_M| +E|r_M|)^i;|r_M| \geq E|r_M|]
\nonumber\\
&&{}  + E[(|r_M| +E|r_M|)^i; |r_M| <E|r_M|]
\\
&\leq& 2^i [E(|r_M|^i)+ (E|r_M|)^i]
\nonumber\\
&\leq& C_4 (2M)^i i!
\nonumber
\end{eqnarray}
By (\ref{e2.13}), (\ref{e2.9}) follows when $0\leq t\leq1/(4M)$.
\end{pf}

Given two points $u=(u_1, u_2)$ and $v=(v_1, v_2)$ with $u_1\leq v_1$ and
$u_2\leq v_2$, we define
$u\stackrel{1}{\rightarrow}v$ as the event that there exists a
northeast 1-path from $u$ to $v$, and define
the \textit{slope} between them by
\[
sl(u, v)={v_2-u_2 \over v_1-u_1}.
\]

%-------------------------------
%f3 ###
\begin{figure}[b]

\includegraphics{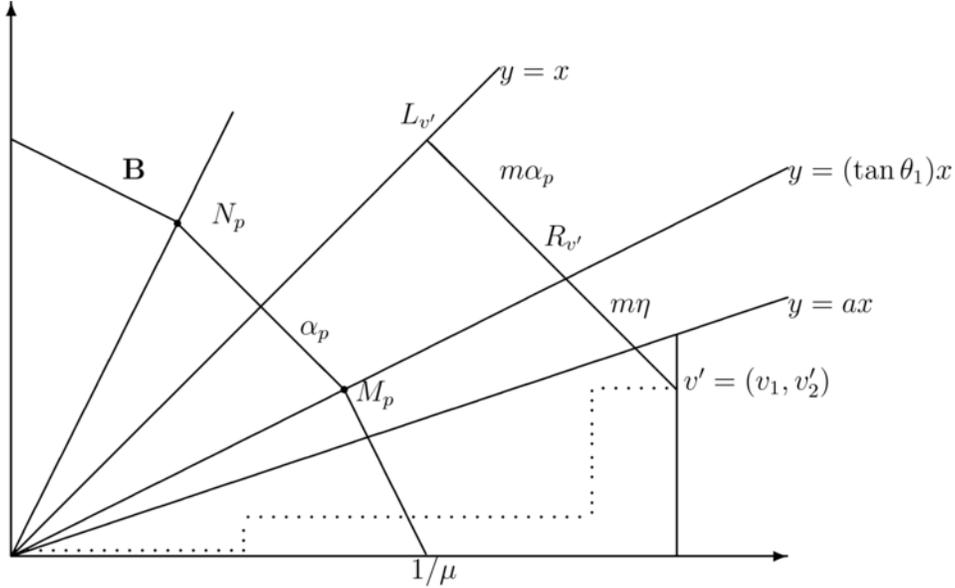}

\caption{The graph shows the relationship between $r_n$ and
$sl((0,0), v')$.
The dotted line is $\gamma_{(0,0), v'}$, the lowest northeast 1-path.
Two lines, $y=x$ and the line perpendicular to $y=x$ passing
through~$v'$, intersect at $L_{v'}$.
The length of the segment from the origin to $L_{v'}$ is $m$. By our
construction,
$m\geq v_1/\sqrt{2}$. We draw a line, $y=(\tan\theta_1) x$, passing
through the origin and
$M_p =({1\over2}+{\alpha_p\over\sqrt{2}}, {1\over2}-{\alpha_p\over \sqrt{2}})$.
The line passes through the above
perpendicular line at $R_{v'}$. By Marchand's result (see Figure~\protect\ref{f1}),
the length of the segment from $L_{v'}$ to $R_{v'}$ is $m \alpha_p$.
We also draw a line $y=ax$ for some $a< \tan(\theta_1)$ defined in
Lemma \protect\ref{lem3}.
By our assumption, $\gamma_{(0,0), v'}$ is below the line $y=ax$.
We denote by $r_m'$ the length of the segment from
$L_{v'}$ to $v'$. Then $r_m'\geq m \alpha_p+ m \eta$.
If we rotate our lattice counterclockwise by $45^\circ$ and extend each
edge by a factor of $\sqrt{2}$,
then $m$ becomes $\sqrt{2}m$ and $ r_{\sqrt{2}m}=r_m'$ for the
right-most edge $r_n$.\label{f3}}
\end{figure}
%%-------------------------------

With these definitions, we show the following lemma.
\begin{lem}\label{lem3}
For $0<a< \tan(\theta_1)$,
if (\ref{e1.10}) holds, then there exist $C_i=C_i(F,a)$ for $i=1,2$ such that
for all $u,v\in\mathbf{Z}^2$,
\[
P\bigl(u\stackrel{1}{\rightarrow}v\mbox{ with }sl(u,v)\leq a\bigr)
\leq C_1\exp\bigl(-C_2 (v_1-u_1)\bigr) .
\]
\end{lem}
\begin{pf}
Suppose that there exists a northeast 1-path from $u$ to
$v$ with $sl(u,v) \leq a$
for some $a < \tan(\theta_1)$. Since we need to use the estimate of
Lemma \ref{lem2}, we introduce the lowest northeast $1$-path.
Let $\gamma_{u,v'}$ be the lowest northeast $1$-path from $u$ to $\{x=v_1\}$,
where the lowest $1$-path means
that all such northeast $1$-paths have no vertex below $\gamma_{u,v'}$.
Let the  last vertex of $\gamma_{u,v'}$
be $v'=(v_1, v_2')$. To show Lemma \ref{lem3}, we need to show
%
%e2.14 ###
\begin{eqnarray}\label{e2.14}
P\bigl(u\stackrel{1}{\rightarrow}v\mbox{ with }sl(u,v)\leq a\bigr)
&\leq & P\bigl(\exists \gamma_{u,v'}, sl(u,v')\leq a\bigr)
\nonumber\\[-8pt]
\\[-8pt]
&\leq & C_1\exp\bigl(-C_2 (v_1-u_1)\bigr).
\nonumber
\end{eqnarray}
The first inequality in (\ref{e2.14}) follows from the definition of the lowest
northeast $1$-path, so we only show the second
inequality in (\ref{e2.14}).
By translation invariance, we may assume that $u=(0,0)$.
Figure~\ref{f3} shows the relationship between $r_n$ and
$sl((0,0), v')$.
By using Figure~\ref{f3}, if $\{\exists \gamma_{u,v'}\mbox{ with }sl((0,0),v')\leq a\}$
for some \mbox{$a< \tan(\theta_1)$}, then there exists $\eta=\eta(a)>0$
such that, for some $m \geq v_1/\sqrt{2}$,
%
%e2.15 ###
\begin{equation}\label{e2.15}
r_m' \geq m (\alpha_p +\eta)\Rightarrow r_n\geq n (\alpha_p +\eta)
\qquad \mbox{for } n=\sqrt{2}m.
\end{equation}
Here, if $n$ is not an integer, we may define $r_n$ as
$r_{\lfloor n\rfloor}$. This implies that
%
%e2.16 ###
\begin{equation}\label{e2.16}
P\bigl(\exists \gamma_{u,v'}\mbox{ with } sl((0,0),v')\leq a\bigr)
\leq\sum_{n \geq v_1/2} P\bigl(r_n \geq n(\alpha_p+\eta)\bigr).
\end{equation}
By (\ref{e2.16}) and Lemma \ref{lem2}, Lemma \ref{lem3} follows.
\end{pf}

Now we show the following two lemmas in order to explore the passage
times in different directions.
\setcounter{equation}{18}%
\begin{lem}[(Marchand \cite{mar02})]\label{lem4}
Under (\ref{e1.10}),
%
%e2.19 ###
\begin{equation}\label{e2.19}
1/\mu> 1/2+\alpha_p/ \sqrt{2}.
\end{equation}
\end{lem}

Recall that the two polar coordinates, for the flat edge on the shape
in the
first quadrant, are denoted by $(\sqrt{1/2+\alpha_p^2}, \theta_i)$ for $i=1,2$,
where
\[
\theta_i= \operatorname{arc}\tan
\biggl( { 1/2\mp\alpha_p/\sqrt{2} \over 1/2\pm \alpha_p/\sqrt{2} } \biggr).
\]
\begin{lem}\label{lem5}
If $F$ satisfies (\ref{e1.10}), there exists $\eta=\eta (F)>0$ such
that
\[
2\geq1+\tan(\theta_1)={\mu+\eta}.
\]
\end{lem}
\begin{pf}
Since
%
%e2.20 ###
\begin{equation}\label{e2.20}
\tan(\theta_1)= {1/2-\alpha_p/\sqrt{2} \over 1/2+\alpha_p/\sqrt{2}},
\end{equation}
then
\[
1+\tan(\theta_1) \leq2.
\]
By (\ref{e2.20}),
%
%e2.21 ###
\begin{equation}\label{e2.21}
\bigl(1+ \tan(\theta_1)\bigr)\bigl(1/2+\alpha_p/\sqrt{2}\bigr)=1.
\end{equation}
By Lemma \ref{lem4}, we take $\eta>0$ such that
%
%e2.22 ###
\begin{equation}\label{e2.22}
1/2+\alpha_p/\sqrt{2}={1\over\mu+\eta},
\end{equation}
so Lemma \ref{lem5} follows from (\ref{e2.22}).
\end{pf}

Now we will introduce two lemmas regarding the rate of convergence of point-
point and point-line passage times.
Kesten \cite{kes93} proved that if $F$ satisfies (\ref{e1.5}) and $F(0) < p_c$,
there exist $C_i= C_i(F)$ for $i=1,2$ such that for all $0< n$ and all
$0< x \leq C_1n$,
%
%e2.23 ###
\begin{equation}\label{e2.23}
P(|\theta_{0,n}-E\theta_{0,n}|\geq x n^{1/2})\leq C_1 \exp(-C_2 x)
\end{equation}
for $\theta=a,b$.
Alexander \cite{ale93} used (\ref{e2.23})
to show that there exist $C= C(F)$ such that for all $0< n$
%
%e2.24 ###
\begin{equation}\label{e2.24}
n\mu\leq Ea_{0,n} \leq n \mu+ C\sqrt{n} \log n.
\end{equation}
If we combine (\ref{e2.23}) and (\ref{e2.24})
together, we have the following lemma.
\begin{lem}[(Alexander \cite{ale93} and Kesten \cite{kes93})]\label{lem6}
If $F$ satisfies (\ref{e1.5}) and $F(0) < p_c$,
there exist $C_i= C_i(F)$ for $i=1,2$ such that for all $0< n$ and all
$n^{0.01}\leq x \leq n^{0.99}$,
\[
P(|a_{0,n}-\mu n|\geq x n^{1/2})\leq C_1 \exp(-C_2 x).
\]
\end{lem}
\begin{lem}\label{lem7}
If $F$ satisfies (\ref{e1.5}) and $F(0)< p_c$,
there exist $C_i= C_i(F)$ for $i=1,2$ such that for all $0< n$ and
$n^{0.01}\leq x\leq n^{0.99}$,
\[
P(|b_{0,n}-\mu n|\geq x n^{1/2})\leq C_1 \exp(-C_2 x).
\]
\end{lem}
\begin{pf}
Lemma 7 was proved in Zhang \cite{zha05}, but the paper was not
published, so here
we reprove it. We may select an optimal path, denoted by $\gamma_n^b$, for
$b_{0,n}$ in a unique way. Then $\gamma^b_n\cap\{x=n\}$
contains only one vertex denoted by $(n, h_n(\gamma_n^b))$.
Smythe and Wierman \cite{smyWie78} proved in their Corollary 8.16 that
\[
\limsup_{n} {h_n(\gamma_n^b)\over n} \leq1 \qquad \mbox{a.s.}
\]
Proposition 5.8 in Kesten \cite{kes86}
tells that, under $F(0) < p_c$, there exist
positive numbers $C_i=C_i(F, \delta)$ for $i=1,2,3$ such that
%
%e2.25 ###
\begin{eqnarray}\label{e2.25}
&& P\bigl(\exists\mbox{ a path $\gamma$ from the origin with $|\gamma| \geq n$ and
$T(\gamma) \leq C_1 n$}\bigr)
\nonumber\\[-8pt]
\\[-8pt]
&&\qquad \leq C_2 \exp(-C_3 n).
\nonumber
\end{eqnarray}
Note that if $h_n(\gamma_n^b) \geq Mn$ for some $M>0$, then
$|\gamma_n^b|\geq Mn$.
It also follows from a large deviation estimate
(see Kesten \cite{kes86}) that there exist
$C_i=C_i(F)$ for $i=3,4$ such that
%
%e2.26 ###
\begin{equation}\label{e2.26}
P\bigl(T(\gamma_n^b) \geq2n \mu\bigr) \leq C_3 \exp(-C_4 n).
\end{equation}
With these observations,
there exist $M=M(F)$ and $C_i(F, M)=C_i$ for $i=5,6$ such that
%
%e2.27 ###
\begin{eqnarray}\label{e2.27}
\qquad
P\bigl(h_n(\gamma_n^b) \geq Mn \bigr)
&\leq & P\bigl( |\gamma_n^b| \geq Mn, T(\gamma_n^b)\leq 2\mu n\bigr)
+ C_3\exp(-C_4 n)
\nonumber\\[-8pt]
\\[-8pt]
&\leq & C_5 \exp(-C_6 n).
\nonumber
\end{eqnarray}
With (\ref{e2.27}),
%
%e2.29 ###
%e2.28 ###
\begin{eqnarray}
&& P(b_{0,n} \leq n\mu-x n^{1/2})
\nonumber\\[-8pt]\label{e2.28}
\\[-8pt]
&& \qquad \leq P\bigl(b_{0,n}
\leq   n\mu-x n^{1/2}, h_n(\gamma_n^b) \leq Mn\bigr)+C_1 \exp(-C_2 n).
\nonumber\\
&& P\bigl(b_{0,n}\leq n \mu-x n^{1/2}, h_n(\gamma_n^b)\leq Mn\bigr)
\nonumber\\[-8pt]\label{e2.29}
\\[-8pt]
&&\qquad \leq  \sum_{i=-Mn}^{Mn}
P\bigl(b_{0,n}\leq n \mu-x n^{1/2}, h_n(\gamma_n^b)=i\bigr).
\nonumber
\end{eqnarray}
From (\ref{e2.28}) and (\ref{e2.29}), there exists $\bar{i}$ such that
%
%e2.30 ###
\begin{eqnarray}\label{e2.30}
\qquad \quad
P(b_{0,n}\leq n \mu-x n^{1/2})
&\leq & 2(M+1)n P\bigl(b_{0,n}\leq n \mu-x n^{1/2}, h_n(\gamma_n^b)=\bar{i}\bigr)
\nonumber\\[-8pt]
\\[-8pt]
&&{} + C_1 \exp(-C_2 n).
\nonumber
\end{eqnarray}
Let $\bar{b}_{0,n}$ be the passage time from $(2n, 0)$ to the line $\{
x=n\}$.
We also select an optimal path $\bar{\gamma}_n^b $ for $\bar{b}_{0,n}$
in a
unique way and denote
\[
(n, \bar{h}_n(\bar{\gamma}_n^b))= \bar{\gamma}_n ^b \cap\{x=n\}.
\]
If
\[
\{b_{0,n}\leq n \mu-x n^{1/2}, h_n(\gamma_n^b)=\bar{i},
\bar{b}_{0,n}\leq
n\mu-x n^{1/2}, \bar{h}_n(\bar{\gamma}_n^b)=\bar{i}\},
\]
then
\[
a_{0,2n} \leq2 n\mu-2x n^{1/2}=2 n\mu-\sqrt{2}x (2n)^{1/2}.
\]
Note that $\{b_{0,n}\leq n \mu-x n^{1/2}, h_n(\gamma_n^b)=\bar{i}\}$ and
$\{\bar{b}_{0,n}\leq n\mu-x n^{1/2}, \bar{h}_n(\bar{\gamma}_n^b)=\bar{i}\}$ only
depend on the configurations of the edges in $-\infty< x < n$
and $n< x < \infty$, respectively,
so the two events are independent and have the same probability. By
Lemma \ref{lem6},
%
%e2.31 ###
\begin{eqnarray}\label{e2.31}
\qquad\quad
&& P(b_{0,n}\leq n\mu-xn^{1/2})^2
\nonumber\\
&&\qquad \leq \bigl(2(M+1)\bigr)^2n^2 P\bigl(b_{0,n}\leq n \mu-x n^{1/2},
h_n(\gamma_n^b)=\bar{i}\bigr)^2
\nonumber\\
&&\qquad \quad {} + C_1 \exp(-C_2n)
\nonumber\\
&&\qquad \leq \bigl(2(M+1)\bigr)^2 n^2 P\bigl(b_{0,n}\leq n \mu-x n^{1/2},
h_n(\gamma_n^b)=\bar{i},
\nonumber\\[-8pt]
\\[-8pt]
&&\hspace*{114pt} \bar{b}_{0,n}\leq n\mu-x n^{1/2},
\bar{h}_n(\bar{\gamma}_n^b)=\bar{i} \bigr)
\nonumber\\
&&\qquad \quad {} +C_1\exp(-C_2n)
\nonumber\\
&&\qquad \leq \bigl(2(M+1)\bigr)^2n^2
P\bigl( a_{0,2n} \leq2 n\mu-\sqrt{2}x (2n)^{1/2}\bigr)
+ C_1\exp(-C_2 n)
\nonumber\\
&&\qquad \leq C_3 \exp(-C_4 x).
\nonumber
\end{eqnarray}
On the other hand, note that $b_{0,n} \leq a_{0,n}$,
so Lemma \ref{lem7} follows from (\ref{e2.31}).
\end{pf}

%s3 ###
\section[Proof of Theorem 1]{Proof of Theorem \protect\ref{thm1}}\label{s3}

In this section, we show Theorem \ref{thm1}.
For the optimal path $\gamma_n$ for $a_{0,n}$ from the origin to $(n,0)$,
we denote by $\gamma_n'(m)$ the piece of $\gamma_n$ from the origin to first
meet the line $\{x=m\}$ for $n^{2/3}\leq m\leq n/2$.
Suppose that $\gamma_n'(m)\cap\{x=m\}= v_n(m)$. The path $\gamma_n$
then goes
from
$v_n(m)$ to $(n,0)$. We denote the last piece by $\gamma_n''(m)$. Clearly,
%
%e3.1 ###
\begin{equation}\label{e3.1}
\gamma_n=\gamma_n'(m)\cup\gamma_n''(m)
\quad \mbox{and}\quad
a_{0,n}= T(\gamma_n'(m)) + T(\gamma_n''(m)).
\end{equation}
If we denote by $\bar{b}_{m,n}$ the passage time from $(n,0)$ to the line
$\{x=m\}$, then
\[
T(\gamma_n''(m) )\geq\bar{b}_{m,n}.
\]
Note that $\bar{b}_{m,n}$ has the same distribution as $b_{0,n-m}$, so
by Lemma \ref{lem7},
%
%e3.2 ###
\begin{eqnarray}\label{e3.2}
\qquad
P\bigl(T(\gamma_n''(m))\leq(n-m) \mu- n^{4/7}\bigr)
&\leq & P\bigl(b_{0, n-m} \leq (n-m)\mu- n^{4/7}\bigr)
\nonumber\\[-8pt]
\\[-8pt]
&\leq & C_1 \exp(-C_2n^{1/14}).
\nonumber
\end{eqnarray}
By (\ref{e3.1}) and (\ref{e3.2}),
%
%e3.3 ###
\begin{eqnarray}\label{e3.3}
&& P\bigl(a_{0,n}\leq T(\gamma_n'(m))+(n-m)\mu-n^{4/7}\bigr)
\nonumber\\
&&\qquad = P\bigl(T(\gamma_n''(m))\leq(n-m) \mu- n^{4/7}\bigr)
\\
&&\qquad \leq C_1 \exp(-C_2n^{1/14}).
\nonumber
\end{eqnarray}
By (\ref{e3.3}),
\begin{eqnarray*}
&& P\bigl(T(\gamma_n'(m)) \geq\mu m + 2n^{4/7}\bigr)
\\
&&\qquad \leq P\bigl(T(\gamma_n'(m)) \geq\mu m + 2n^{4/7},
a_{0,n}\geq T(\gamma_n'(m))+(n- m)\mu-n^{4/7}\bigr)
\\
&&\qquad \quad {} +C_1\exp(-C_2 n^{1/14}).
\end{eqnarray*}
Note that $\{T(\gamma_n'(m)) \geq\mu m + 2n^{4/7},a_{0,n}\geq
T(\gamma_n'(m))+(n-m)\mu-n^{4/7}\}$
implies that $a_{0,n}\geq n\mu+n^{4/7}$, so by Lemma \ref{lem6},
%
%e3.4 ###
\begin{eqnarray}\label{e3.4}
\qquad 
P\bigl(T(\gamma_n'(m)) \geq\mu m + 2n^{4/7}\bigr)
&\leq & P(a_{0,n}\geq n\mu +n^{4/7})+C_1\exp(-C_2 n^{1/14})
\nonumber\\[-8pt]
\\[-8pt]
&\leq & C_3 \exp(-C_4 n^{1/14}).
\nonumber
\end{eqnarray}
Now we estimate the length of $\gamma_n'(m)$.
Clearly, $|\gamma_n'(m)| \geq m$.
Note that $m \geq n^{2/3}> \mu^{-1}n^{4/7}$ for large $n$, so
by (\ref{e3.4}) we have
%
%e3.5 ###
\begin{equation}\label{e3.5}
P\bigl(T(\gamma_n'(m)) \geq3\mu m \bigr)\leq C_1\exp(-C_2 n^{1/14}).
\end{equation}
As $t(e) \geq1$ almost surely for all edges, if $N=3\mu$, then on
$T(\gamma_n'(m)) \leq3\mu m$,
%
%e3.6 ###
\begin{equation}\label{e3.6}
|\gamma_n'(m) |\leq N m \qquad \mbox{almost surely}.
\end{equation}
Together with (\ref{e3.5}) and (\ref{e3.6}), we have for all $n$ and $m \geq n^{2/3}$,
%
%e3.7 ###
\begin{equation}\label{e3.7}
P\bigl(|\gamma_n'(m) |\geq N m\bigr) \leq C_1\exp(-C_2 n^{1/14}).
\end{equation}

Now we use the method of renormalization in Kesten and Zhang \cite{kesZha90}.
We define, for integer $M$ and $w=(w_1, w_2)\in\mathbf{Z}^2$, the squares
and the vertical strips by
\[
B_M(w) = [Mw_1, Mw_1+M)\times[Mw_2, Mw_2+M)
\]
and
\[
V_M(w_1)=[Mw_1, Mw_1+M)\times(-\infty, \infty).
\]
We denote the $M$-squares and the $M$-strip by
$\{B_M(w)\dvtx w\in\mathbf{Z}^2\}$ and $\{V_M(w_1)\dvtx\break w_1 \in\mathbf{Z}\}$,
respectively. For the optimal path $\gamma_n$, we denote
a fattened ${\gamma}_n'(m,M)$ by
\[
{\gamma}_n'(m,M)=\{B_M(w)\dvtx B_M(w)\cap\gamma_n'(m)\neq\varnothing\}.
\]
By our definition,
%
%e3.8 ###
\begin{equation}\label{e3.8}
|{\gamma}_n'(m,M)|\geq{|\gamma_n'(m)|\over M^2}\geq{m\over M^2}.
\end{equation}
For each $M$-square $B_M(w)$, there are eight $M$-square neighbors. We
say they are adjacent to $B_M(w)$.
Since $\gamma_n(m)'$ is connected, ${\gamma}_n'(m,M)$ has to be connected
through the square connections.

Note that if there are much fewer
vertices of $[S_M(\gamma_n)\cup D(\gamma_n)]$ than other vertices in
$\gamma_n'(m)$,
then there are also fewer strips that contain a vertex in
$[S_M(\gamma _n)\cup D(\gamma_n)]$ than the other strips.
We say a strip $V_M(w_1)$ is \textit{bad} if
\[
B(m)\cap[S_M(\gamma_n)\cup D(\gamma_n)]\cap V_M(w_1)\neq\varnothing.
\]
For each bad strip $V_M(w_1)$, we also say two neighbor strips to its
left and
two neighbor strips to its right
are bad. Otherwise, we say a strip is \textit{good}.

\begin{figure}

\includegraphics{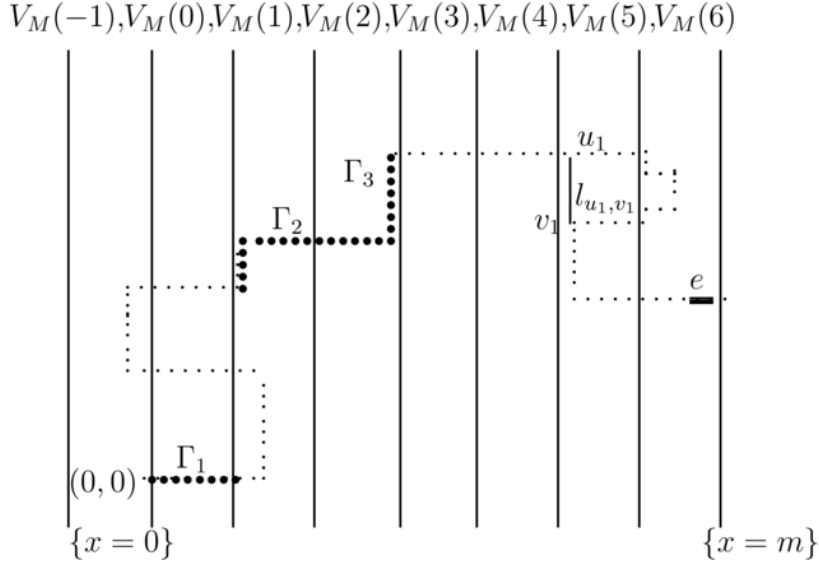}

\caption{The graph shows the path $\gamma_n'(m)$ in the strips,
where the dotted line is $\gamma_n'(m)$.
$\Gamma_m$ is the path from the origin to $u_1$ along the dotted line
from $u_1$ to $v_1$ along
$l_{u_1, v_1}$ and then from $v_1$ along the dotted line to $\{x=m\}$.
$\Gamma_m$ crosses three good
strips: $V_M(0)$ with the big dotted line $\Gamma_1=\Gamma_m(v_M(0,0),
v_M(0,1))$; $V_M(1)$ with
the big dotted line $\Gamma_2=\Gamma_m(v_M(1,0),v_M(1, 1))$; and
$V_M(2)$ with the big dotted line $\Gamma_3=\Gamma_m(v_M(2,0),v_M(2, 1))$.
$V_M(0)$ is a good-short-flat strip since the big dotted line $\Gamma_1$ is
flat. $V_M(1)$ is a good-short-nonflat strip since the slope of its two terminal
points in $\Gamma_2$ is larger than
$\tan(\theta_1)-\delta_1$. $V_M(2)$ is a good-long strip since the
number of vertices in $\Gamma_3$ is larger than $2M$.
$V_M(5)$ is bad because it contains
$M$-broken bridge $l_{u_1,v_1}$ for $\gamma_n$.
$V_M(3)$ and $V_M(4)$ would be good-short-flat strips, however, they
are bad since they are next to
$V_M(5)$. $V_M(6)$ is bad since $e$ is a $1^+$-edge, and it is also
next to bad strip $V_M(5)$. After eliminating the bad strips,
$\mathbf{G} = V_M(0)\cup V_M(1)\cup V_M(2)$.\label{f4}}
\end{figure}

We eliminate all bad strips from $\mathbf{Z}^2$ (see Figure~\ref{f4})
and denote the remaining vertices by
${\bf G}$. Recall our definitions of $\{l_{u_i, v_i}\}$ and $\gamma(u_i, v_i)$
for $i=1,\ldots,\tau$, in Section \ref{s1}. We may define
$\Gamma_m$ (see Figure~\ref{f4})
as the path from the origin that goes along $\gamma_n$, meets $u_i$,
then goes along $l_{u_i, v_i}$
from $u_i$ to $v_i$ (not along $\gamma(u_i, v_i)$ in $\gamma_n$), and then
goes along $\gamma_n$ from $v_i$ to $u_{i+1}$,
until it meets $\{x=m\}$ for $i=1,\ldots,\tau$.

%-----------------------------------------------
%f4 ###
%-----------------------------------------------

For a strip with
\[
V_M(w_1)\subset\{0< x< m\},
\]
$\Gamma_m$ will cross the strip $V_M(w_1)$.
If we go along $\Gamma_m$ to cross through $V_M(w_1)$, we will meet vertex
$v_M(w_1,1)$ at the left boundary of $V_M(w_1)$, and then go along
$\Gamma_m$
using the vertices inside $V_M(w_1)$ to meet vertex $v_M(w_1, 2)$ at
the right
boundary of
$V_M(w_1)$. Note that $\Gamma_m$ may cross through $V_M(w_1)$ back and forth
many times (see Figure~\ref{f4}). We select one of them
(see Figure~\ref{f4} for $\Gamma_1$) and denote this subpath of $\Gamma_m$ from
$v_M(w_1,1)$ to $v_M(w_1,2)$ by $\Gamma_m(v_M(w_1,1), v_M(w_1,2))$.
If $V_M(w_1)\subset{\bf G}$, $\Gamma_m(v_M(w_1,1), v_M(w_1,2))$ cannot
contain a vertex of
$l_{u_i, v_i}$ for $i=1,\ldots, \tau$, so
\[
\Gamma_m(v_M(w_1,1), v_M(w_1,2)) \subset\gamma_n(m)'.
\]
For each such path $\Gamma_m(v_M(w_1,1), v_M(w_1,2))$, the path cannot
contain a
vertex $v$ such that $v\in [S_M(\gamma_n)\cup D(\gamma_n)]$ by our definition.
This tells us that $\Gamma_m(v_M(w_1,1), v_M(w_1,2))$ is a $1$-path.
If
\[
|\Gamma_m(v_M(w_1,1), v_M(w_1,2))|\geq2M,
\]
we say the strip $V_M(w_1)$ is a \textit{good-long strip}.
Otherwise, it is a \textit{good-short strip}.

Now we focus on all good-short strips.
Assume $V_M(w_1)$ is a good-short strip.
Then by (\ref{e1.15}), $\Gamma_m(v_M(w_1,1), v_M(w_1,2))$ has none of its own
$M$-broken bridges. By Lemma~\ref{lem1},
$\Gamma_m(v_M(w_1,1), v_M(w_1,2))$ is either northeast or southeast.
We say that a good-short strip $V_M(w_1)$ is \textit{good-short-flat} if
\[
|sl(v_M(w_1,1), v_M(w_1,2) )|\leq\tan(\theta_1)-\delta_1\mbox,
\]
where $0<\delta_1< \tan(\theta_1)$ is taken such that
%
%e3.9 ###
\begin{equation}\label{e3.9}
\delta_1 = \eta/100
\end{equation}
for the $\eta$ defined in Lemma \ref{lem5}.
In contrast, $V_M(w_1)$ is \textit{good-short-non-flat} if
\[
|sl(v_M(w_1,1), v_M(w_1,2) )|> \tan(\theta_1)-\delta_1.
\]

For a good-short-flat strip $V_M(w_1)$, if $v_M(w_1,1)\in B_M(w)$ for some
$w=(w_1,w_2)$,
then we say the square $B_M(w)$ is a
\textit{good-short-flat square}.
By Lemma~\ref{lem3}, for a fixed $B_M(w)$, there exist positive constants
$\beta_i=\beta_i(F, \delta_1)$ for $i=1,2,3,4$ such that for all $M$,
%
%e3.10 ###
\begin{eqnarray}\label{e3.10}
&& P(B_M(w) \mbox{ is a good-short-flat square})
\nonumber\\[-8pt]
\\[-8pt]
&&\qquad \leq\beta_1 M \exp(-\beta_2 M)\leq\beta_3 \exp(-\beta_4 M).
\nonumber
\end{eqnarray}

Now we denote by $F(\gamma_n'(m))$ the number of good-short-flat
strips. We shall show that
there exist $\delta_2 >0$ with $\delta_2 \leq \eta/100 $ for the
$\eta$ in Lemma~\ref{lem5}, $M=M(\delta_2, \beta_3, \beta_4, N)$, and
$C_i=C_i(F, N, M,\delta_1,\delta_2)$ for
$i=1,2$ such that for all $m$ and $n$ with $n^{2/3} \leq m\leq n/2$,
%
%e3.11 ###
\begin{equation}\label{e3.11}
P\bigl(F(\gamma_n'(m)) \geq\delta_2 m/M\bigr)\leq C_1 \exp(-C_2
n^{1/14}).
\end{equation}

To show (\ref{e3.11}), we need to introduce a few basic methods to account for
connected squares. By (\ref{e3.7}),
%
%e3.12 ###
\begin{eqnarray}\label{e3.12}
\qquad
&& P\bigl(F(\gamma_n'(m)) \geq\delta_2 m/M\bigr)
\nonumber\\[-8pt]
\\[-8pt]
&&\qquad \leq P\bigl(F(\gamma_n'(m)) \geq\delta_2 m/M,
|\gamma_n'(m)| \leq Nm\bigr)+C_1 \exp(-C_2 n^{1/14}).
\nonumber
\end{eqnarray}

As we mentioned, ${\gamma}_n'(m,M)$ is connected through horizontal, vertical,
and diagonal squares. Therefore,
there are at most $8^{k}$ choices for all $M$-squares in the path
${\gamma}_n'(m,M)$ if
$|{\gamma}_n'(m,M)|=k$, where $|{\gamma}_n'(m,M)|$ represents the
number of $M$-squares in ${\gamma}_n'(m,M)$.
Let $\bar{B}_M(w)$ be the union of $B_M(w)$ and its eight neighbor $M$-squares.
We call this a $3M$-\textit{square}.
If $B_M(w)\cap\gamma_n'(m)\neq\varnothing$, then $\bar{B}_M(w)$
contains at
least $M$ vertices of $\gamma_n'(m)$
in its interior. We collect all such $3M$-squares $\{ \bar{B}_M(w)\}$ such
that
their center $M$-squares contain at least a vertex of $\gamma_n'(m)$.

We need to decompose $\gamma'_n(m, M)$ into disjoint $3M$-squares.
We select an $M$-strip, from $\{V_M(w_1)\dvtx w_1\in\mathbf{Z}\}$,
with the maximum number of $M$-squares in $\gamma_n'(m, M)$.
We denote by $V_{\max}(1)$ and $i_1$ the strip and the number of $M$-squares
in $V_{\max}(1)\cap\gamma_n'(m, M)$,
respectively. We also denote by
$V^-_{\max}(1)$ the two neighboring $M$-strips on the left side of
$V_{\max}(1)$.
In addition, we denote by $i_1^-$ the number of $M$-squares in
$V_{\max}^-(1)\cap\gamma_n'(m, M)$. Similarly,
$V_{\max}^+(1)$ and $i^+_1$ are denoted by the two neighboring
$M$-strips on the right side of $V_{\max}(1)$, and the number
of the $M$-squares in $V_{\max}^+(1)\cap\gamma_n'(m, M)$, respectively.
By this construction,
\[
i_1^+\leq2i_1 \quad \mbox{and} \quad i_1^-\leq2i_1.
\]
We consider the top square,
denoted by $B_M(u_1)\in{\gamma}_n'(m,M)$, in $M$-strip $V_{\max}(1)$.
Clearly, $\bar{B}_M(u_1)$ contains
at least $M$ vertices in $\gamma_n'(m)$. Now we consider the second top
$M$-square in this strip, denoted by
$B_M(u_2)\subset{\gamma}_n'(n,M)$, such that
$\bar{B}_M(u_2)$ and $\bar{B}_M(u_1)$ have no common vertices. We
continue in
this way to find all the disjoint $3M$-squares in this strip
such that their center $M$-squares contain at least a vertex of $\gamma_n'(m)$.
Note that it is possible that only one $3M$-square exists in this strip.

Next we select an $M$-strip, from $\{V_M(w_1)\dvtx
w_1\in\mathbf{Z}\}\setminus \{V_{\max}(1)\cup V_{\max}^\pm(1)\}$,
with the maximum number of $M$-squares in $\gamma_n'(m, M)$.
We denote by $V_{\max}(2)$ and $i_2$ the strip and the number of $M$-squares
in $V_{\max}(2)\cap\gamma_n'(m, M)$,
respectively.
We also select the two neighboring $M$-strips on the left side $V_{\max}(2)$.
Note that these two strips might overlap
the strips of $V_{\max}^-(1)\cup V_{\max}^+(1)$. If so, we eliminate the
overlapped strips from these two strips.
We denote the selected strips by
$V^-_{\max}(2)$. In addition, we denote by $i_2^-$ the number of $M$-squares
in $V_{\max}^-(2)\cap\gamma_n'(m, M)$. Note that $V^-_{\max}(2)$ might be empty
after the overlapped strips are eliminated.
If it is empty, then $i_2^-=0$. Similarly, we select
the two neighboring $M$-strips on the right side of $V_{\max}(2)$. After
eliminating the overlapped strips of
$V_{\max}^-(1)\cup V_{\max}^+(1)$ from these two strips, let
$V_{\max}^+(2)$ and $i^+_2$ denote the selected strips and the number
of $M$-squares in $V_{\max}^+(2)\cap\gamma_n'(m, M)$, respectively.
By this construction,
\[
i_2^+\leq2i_2 \quad \mbox{and}\quad  i_2^-\leq2i_2.
\]
With $V_{\max}(2)$, we select $3M$-squares in the same way that we
selected for $V_{\max}(1)$.
Note that the selected $3M$-squares in $V_{\max}(1)$ and $V_{\max}(2)$
are disjoint.

We then continue this process to select from the third to the last $M$-strip
to find all the disjoint $3M$-squares such that their center $M$-squares contain
at least a vertex of $\gamma_n'(m)$. By these selections,
we have
%
%e3.13 ###
\begin{equation}\label{e3.13}
i_j^+\leq2i_j \quad \mbox{and}\quad
i_j^-\leq2i_j \qquad \mbox{for all }j.
\end{equation}
By (\ref{e3.13}),
\[
k= \sum_{j} (i_j+i_j^++i_j^-)\leq\sum_{j} 5 i_j.
\]
Thus,
\[
k/5\leq\sum_{j}i_j.
\]
By our construction,
each strip $V_{\max}(j)$ contains at least $\lceil i_j/3\rceil$
of these disjoint $3M$-squares.
With these observations, if $|\gamma'(m,M)|=k$ for large $k$, there are
at least $k/15$ disjoint $3M$-squares
such that their center $M$-squares contain
at least a vertex of $\gamma_n'(m)$.

If there are at most $Nm$ vertices in $\gamma_n'(m)$, then
\[
(Mk/15)\leq Nm.
\]
On the other hand, note that $\gamma_n'(m)$ is a path from $(0,0)$ to
$\{x=m\}$,
so
it at least crosses $m/M$ strips. This implies that
\[
|{\gamma}_n'(m,M)|\geq m/M.
\]
With these observations, on $|\gamma_n'(m) |\leq Nm$,
%
%e3.14 ###
\begin{equation}\label{e3.14}
m/M\leq|{\gamma}_n'(m,M)|=k \leq15 Nm/M.
\end{equation}

By (\ref{e3.14}),
%
%e3.15 ###
\begin{eqnarray}\label{e3.15}
\qquad \quad
&& P\bigl(F(\gamma_n'(m)) \geq\delta_2 m/M, |\gamma_n'(m)|\leq N m\bigr)
\nonumber\\
&&\qquad = \sum_{15Nm/M \geq k\geq m/M}
P\bigl(F(\gamma_n'(m)) \geq\delta_2 m/M,
\nonumber\\
&& \hspace*{117pt}
|\gamma_n'(m)|\leq N m,|{\gamma}_n'(m,M)|=k \bigr)
\\
&&\qquad \leq \sum_{15Nm/M \geq k\geq m/M}8^k
P\bigl(F(\gamma_n'(m)) \geq\delta_2 m/M,
\nonumber\\
&&\hspace*{128pt}
|\gamma_n'(m)|\leq N m,{\gamma}_n'(m,M)=\Gamma, |\Gamma|=k\bigr),
\nonumber
\end{eqnarray}
where $\Gamma$ is a fixed connected $M$-squares and $|\Gamma|=k$ means that
$\Gamma$ contains $k$ squares.
On the event $\{F(\gamma_n'(m)) \geq\delta_2 m/M, |\gamma_n'(m)|\leq N m,
{\gamma}_n'(m,M)=\Gamma, |\Gamma|=k)\}$, there exists a connected set
$\Gamma$ of $k$ $M$-squares, with the bounds in (\ref{e3.14}) for~$k$, and that
contains at least $\delta_2 m/M$ good-short-flat squares.
If there are more than $\delta_2 m/M$ of such good-short-flat strips,
note that $v_M(w_1, 1)$ has to stay in the
left boundary of an $M$-square $B_M(w)\subset{\gamma}_n'(m,M)$, so we select
all such $\delta_2m/M$ good-short-flat squares from $\Gamma$ to have at most
%
%e3.16 ###
\begin{equation}\label{e3.16}
\pmatrix{ 15Nm/M \cr
\delta_2 m/M}\leq2^{15Nm/M}\mbox{ choices}.
\end{equation}
Therefore, by (\ref{e3.10}), (\ref{e3.15}) and (\ref{e3.16}),
%
%e3.17 ###
\begin{eqnarray}\label{e3.17}
&& P\bigl(F(\gamma_n'(m)) \geq\delta_2 m/M, |\gamma_n|\leq N m\bigr)
\nonumber\\[-8pt]
\\[-8pt]
&&\qquad
\leq(15 N m) 8^{15Nm/M} 2^{15Nm/M} \beta_3^{\delta_2m/M} \exp(-\beta_4 m \delta_2 ).
\nonumber
\end{eqnarray}
Thus, for $\eta$ in Lemma 5, $0< \delta_2 \leq\eta/100$, and $N$ in
(\ref{e3.7}), we
select $M=M(\eta, \delta_2, \beta_3, \beta_4, N)$ large such that
(\ref{e3.11}) follows from (\ref{e3.12}) and (\ref{e3.17}).

Now we show that (\ref{e3.11}) implies Theorem \ref{thm1}.
Note that $\gamma_n'(m)$ crosses out from $\{x=0\}$ to $\{x=m\}$.
Without loss of generality, we may assume that $\gamma_n$ first meets
the left or right boundary of $B(m)$.
If not, we can always work on the horizontal strips rather than the vertical
strips by using the same argument.
There are at least $m/M$ strips that have a common vertex with $\gamma_n'(m)$ on
$B(m)$. On
\[
|B(m) \cap[S_M(\gamma_n)\cup D(\gamma_n)]|\leq\delta m
\]
for $\delta= \delta_2/M$ and $\delta_2$ defined above, there are at most
$(5\delta_2m)/M$ bad strips.
After eliminating these bad strips, we have at least $(1-5\delta_2)m/M$ good
strips left.
Note that $\gamma_n'(m)$ contains at least $m$ horizontal edges and
each of
them at least costs time one, so
\[
T(\mbox{horizontal edges of } \gamma_n'(m)) \geq m.
\]
Now we account for the vertical edges in $\gamma_n'(m)$. In fact, under
$F(\gamma_n'(m))\leq\delta_2 m/M$ for a small
$\delta_2>0$, the number of good-short-flat strips is less than $\delta_2 m/M$.
The total number of good-long
and good-short-nonflat strips is more than
$m/M-5\delta_2 m/M -\delta_2 m/M=(1-6\delta_2)m/M$.
For each good-long or good-short-nonflat strip,
by the definition of $\delta_1$ in (\ref{e3.9}), $\Gamma(v_M(u,1), v_M(u,2))$
has to
contain at least
\[
\bigl(\tan(\theta_1) -\delta_1\bigr)M-1
\geq \bigl(\tan(\theta_1) - 0.02\eta\bigr)M
\]
vertical edges, where $-1$ is a possible error made by a
noninteger number $(\tan(\theta_1)-\delta_1)M$.
Therefore, on $F(\gamma_n'(m)) \leq\delta_2 m/M$,
there are at least
\[
[m(1-6\delta_2)/M]\bigl[\bigl(\tan(\theta_1)-0.02\eta\bigr)M\bigr]
= m(1-6\delta_2)\bigl(\tan (\theta_1)- 0.02\eta\bigr)
\]
total vertical edges. The total time of vertical edges costs at least
$(\tan(\theta_1)-0.1\eta)m$. Together with horizontal edges, by
Lemma~\ref{lem5}, on $F(\gamma_n'(m)) \leq\delta_2 m/M$,
%
%e3.18 ###
\begin{equation}\label{e3.18}
T(\gamma_n'(m))\geq m+[\tan(\theta_1)-0.1\eta]m\geq m(\mu+\eta/2).
\end{equation}
Therefore, for $n$ and $m$ with $m \geq n^{2/3} > 4\eta^{-1} n^{4/7}$,
by (\ref{e3.4}), (\ref{e3.11}) and (\ref{e3.18}), there exist
$C_i=C_i(F, \delta,\eta)$ for $i=1,2,3,4$ such that
\begin{eqnarray*}
&& P\bigl(| B(m)\cap[S_M(\gamma_n)\cup D(\gamma_n)]| \leq\delta m\bigr)
\\
&&\qquad \leq2 P\bigl(| B(m)\cap[S_M(\gamma_n)\cup D(\gamma_n)]| \leq\delta
m,F(\gamma_n'(m)) \leq\delta_2 m/M\bigr)
\\
&&\qquad \quad {} +C_1\exp(-C_2 n^{1/14})
\\
&&\qquad \leq 2P\bigl(T(\gamma_n'(m))\geq m(\mu+\eta/2)\bigr)
+ C_1\exp(-C_2 n^{1/14})
\\
&&\qquad \leq C_3 \exp(-C_4 n^{1/14}),
\end{eqnarray*}
where factor 2 above is the result of the assumption that
$\gamma_n'(m)$ meets
the left or right boundary of $B(m)$ first.
So Theorem \ref{thm1} follows.

%s4 ###
\section[Corollaries of Theorem 1]{Corollaries of Theorem \protect\ref{thm1}}\label{s4}
\setcounter{equation}{-1}

In this section, we need to generalize Theorem~\ref{thm1}.
Let $\delta_3$ be a number such that
%
%e4.1 ###
\begin{equation}\label{e4.0}
\delta_3^{\delta/16} (48M)^{N}=1/2,
\end{equation}
where $\delta$, $N$ and $M$ are the numbers selected in Theorem \ref{thm1},
(\ref{e3.6}) and (\ref{e3.11}) in Section \ref{s3}.
Since $F$ is a right continuous function and $t(e)$ is not a constant,
we may select $z>1$ such that
%
%e4.2 ###
\begin{equation}\label{e4.1}
P\bigl(1< t(e) \leq z\bigr)=F(z)- F(1) \leq\delta_3,
\end{equation}
where $\delta_3$ is the number in (\ref{e4.0}). We say $e$ is a $z^+$-edge if
$t(e) >z$, where $z>1$ is the number in (\ref{e4.1}).
For the optimal path $\gamma_n$, we denote by
$D(z,\gamma_n)$ all the vertices in $\gamma_n$ that are adjacent to $z^+$-edges
on $\gamma_n$. We also let
$S_M(z,\gamma_n)$ be the set of vertices in $\gamma_n$ that are
adjacent to
$M$-broken bridges $\{l_{u_i, v_i}\}$ of
$\gamma_n$ and, in addition, each broken bridge $l_{u_i, v_i}$ contains
at least
one $z^+$-edge.
With these definitions, we have the following corollary.
\begin{cor}\label{cor1}
If $F$ satisfies (\ref{e1.5}) and (\ref{e1.10}), there exist
$C_i=C_i(F, M, N,\break z,\delta_3,\delta)$ for $i=1,2,3$, and $M$ in
(\ref{e3.11}), $N$ in
(\ref{e3.7}), and $\delta$ in Theorem \ref{thm1}
and $\delta_5$ in~(\ref{e4.0}), and $z$ in (\ref{e4.1}) such
that
\[
P\bigl(|B(m)\cap[D(z,\gamma_n)\cup S_M(z,\gamma_n)]|\leq\delta m/4\bigr)
\leq C_1\exp(- C_2n^{1/14})
\]
for all $m$ with $n/2 \geq m \geq n^{2/3}$
and
\[
E\bigl(|B(m)\cap[D(z,\gamma_n)\cup S_M(z,\gamma_n)]|\bigr) \geq C_3m.
\]
\end{cor}
\begin{pf}
By Theorem \ref{thm1},
%
%e4.3 ###
\begin{eqnarray}\label{e4.2}
&& P\bigl(|B(m)\cap[D(z,\gamma_n)\cup S_M(z,\gamma_n)]|\leq\delta m/4 \bigr)
\nonumber\\
&&\qquad \leq P\bigl(|B(m)\cap[D(z,\gamma_n)\cup S_M(z,\gamma_n)]|\leq\delta m/4,
\nonumber\\[-8pt]
\\[-8pt]
&&\hspace*{79pt}
|B(m)\cap[D(\gamma_n)\cup S_M(\gamma_n)]|\geq\delta m\bigr)
\nonumber\\
&&\qquad \quad {} +C_1\exp(-C_2n^{1/14}).
\nonumber
\end{eqnarray}
Note that if $|B(m)\cap D(\gamma_n)\cup S_M(\gamma_n)|\geq\delta m $, then
either
\[
|D(\gamma_n)\cap B(m)|\geq\delta m/2,
\]
or
\[
|S_M(\gamma_n)\cap B(m)|\geq\delta m/2.
\]
We may assume that the first event occurs.
For $v\in D(\gamma_n)\cap B(m)$, $v$ is adjacent to $e$ on $\gamma_n$
with $t(e)>1$. Thus, there are at least half of these
vertices in $D(\gamma_n)\cap B(m)$ such that the edges adjacent to these
vertices cannot take a value larger than $z$ under
\[
|B(m)\cap[D(z,\gamma_n)\cap S_M(z, \gamma_n)]|\leq\delta m/4.
\]
In other words, there are at least $\delta m/4$ vertices in $\gamma_n$
that are adjacent to edges~$\{e\}$ on $\gamma_n$
with $1< t(e) \leq z$. Therefore, there are at least $\delta m/16$
edges in $\gamma_n'(m)$ adjacent to these vertices
with $1< t(e) \leq z$, since each vertex is adjacent to at most four
edges. Recall that $\gamma_n'(m)$, defined in the last section, is the
piece of $\gamma_n$ from the origin to
the line $\{x=m\}$.
To fix our path $\gamma_n'(m)$, on $|\gamma_n'(m) |\leq Nm$, we have at most
$4\cdot3^{Nm}$ choices.
After fixing our path $\gamma_n'(m)$, we fix these edges with
$1< t(e) \leq z$, so we have at most
\[
\sum_{k=1}^{Nm} \pmatrix{Nm \cr k} \leq2^{Nm}\mbox{ choices}.
\]
With these observations, by (\ref{e3.7}), if we take $\delta_3$ satisfying (\ref{e4.0}),
%
%e4.4 ###
\begin{eqnarray}\label{e4.3}
\qquad
&& P\bigl(|B(m)\cap[D(z,\gamma_n)\cup S_M(z,\gamma_n)]|
\leq\delta m/4,|B(m)\cap D(\gamma_n)|\geq\delta m/2\bigr)
\nonumber\\
&&\qquad \leq P\bigl(|B(m)\cap[D(z,\gamma_n)\cup S_M(z,\gamma_n)]|\leq\delta m/4,
\nonumber\\
&& \hspace*{46pt}
|B(m)\cap D(\gamma_n)|\geq\delta m/2, |\gamma_n'(m)|\leq Nm\bigr)
\nonumber\\[-8pt]
\\[-8pt]
&&\qquad \quad {} + C_1 \exp(-C_2 n^{1/14})
\nonumber\\
&&\qquad \leq C_3 3^{Nm} 2^{Nm} \delta_3^{\delta m/16}+C_1\exp(-C_2 n^{1/14})
\nonumber\\
&&\qquad \leq C_4 \exp(-C_5 n^{1/14}).
\nonumber
\end{eqnarray}

Now we assume that the second event occurs:
\[
|S_M(\gamma_n)\cap B(m)|\geq\delta m/2.
\]
By (\ref{e1.18}), each $M$-broken bridge has at least an edge $e$
with $t(e)>1$. If this edge is not a $z^+$-edge,
then we have $1< t(e) \leq z$. Note that each $M$-bridge has at most
$2M$ edges. Note also that
if $u$ is fixed, then there are at most four choices for $l_{u,v}$, so
we use the same estimate
as (\ref{e4.3}) to fix the path $\gamma_n'(m)$, the starting vertices
in~$\gamma_n$ for $M$-broken bridges,
the $M$-broken bridges, and
the edges with $1< t(e) \leq z$ in these $M$-broken bridges, resulting in
\begin{eqnarray*}
&& P\bigl(|B(m)\cap[D(z,\gamma_n)\cup S_M(z,\gamma_n)]|\leq\delta m/4,
\\
&&\hspace*{14pt}
|B(m)\cap S_M(\gamma_n)|\geq\delta m/2, |\gamma_n'(m)|\leq Nm\bigr)
\\
&&\qquad \leq C_3 3^{Nm} 2^{Nm} 4^{Nm} (2M)^{Nm} \delta_3^{\delta m/16}.
\end{eqnarray*}
By (\ref{e4.0}), there exist
$C_i=C_i(F, M, N, \delta_3,\delta, z)$ for $i=1,2$ such that
%
%e4.5 ###
\begin{eqnarray}\label{e4.4}
&& P\bigl(|B(m)\cap[D(z,\gamma_n)\cup S_M(z,\gamma_n)]|\leq\delta m/4,
\nonumber\\
&&\hspace*{14pt}
|B(m)\cap S_M(\gamma_n)|\geq\delta m/2, |\gamma_n'(m)|\leq Nm\bigr)
\nonumber\\[-8pt]
\\[-8pt]
&&\qquad \leq C_1\exp(-C_2 m).
\nonumber
\end{eqnarray}
Together with (\ref{e4.3}) and (\ref{e4.4}),
%
%e4.6 ###
\begin{eqnarray}\label{e4.5}
&& P\bigl(|B(m)\cap[D(z,\gamma_n)\cup S_M(z,\gamma_n)]|\leq\delta m/4,
\nonumber\\
&&\hspace*{47pt}|B(m)\cap [D(\gamma_n)\cup S_M(\gamma_n)]|\geq\delta m\bigr)
\nonumber\\
&&\qquad \leq P\bigl(|B(m)\cap[D(z,\gamma_n)\cup S_M(z,\gamma_n)]|\leq\delta m/4,
\nonumber\\
&&\hspace*{52pt}
|B(m)\cap D(\gamma_n)|\geq\delta m/2, |\gamma_n'(m)|\leq Nm\bigr)
\nonumber\\[-8pt]
\\[-8pt]
&&\qquad \quad {}
+ P\bigl(|B(m)\cap[D(z,\gamma_n)\cup S_M(z,\gamma_n)]|\leq\delta m/4,
\nonumber\\
&&\hspace*{60pt}
|B(m)\cap S_M(\gamma_n)|\geq\delta m/2, |\gamma_n'(m)|\leq Nm\bigr)
\nonumber\\
&&\qquad \quad {} + 2 C_1 \exp(-C_2 n^{1/14})
\nonumber\\
&&\qquad \leq C_4 \exp(-C_5 n^{1/14}).
\nonumber
\end{eqnarray}
Therefore, the probability estimate in Corollary \ref{cor1}
follows from (\ref{e4.5}). With (\ref{e4.5}), we also have
%
%e4.7 ###
\begin{equation}\label{e4.6}
E\bigl(|B(m)\cap[D(z,\gamma_n)\cup S_M(z,\gamma_n)]|\bigr) \geq C_3m .
\end{equation}\upqed
\end{pf}

In Corollary \ref{cor1},
we showed that there are proportionally many vertices $\{v\}$ such that
$v\in
S_M(z,\gamma_n)\cup D(z,\gamma_n)$.
If $u\in S_M(z,\gamma_n)$ for the optimal path $\gamma_n$, then to show
Theorem~\ref{thm2},
we need the $M$-broken bridge $l_{u,v}$ to stay inside a large square.
Let us consider $B_{M}(u)$. The square has the center at $(u_1+M/2, u_2+M/2)$.
Now we construct a larger square
$G(B_M(u))$ with the same center at $(u_1+M/2, u_2+M/2)$ and a side
length of
$7M$. Note that $G(B_M(u))$ contains 49 of these \mbox{$M$-squares} and
$B_M(u)$ is the center $M$-square among these 49 \mbox{$M$-squares}. Here we require
these $G$-squares,
the same as for the $B$-squares, to have lower and left boundaries but
no top and right boundaries.

$G(B_M((0,0)))$ contains 49 of these $M$-squares. We denote them by
$\{B_M(q_1),\break\ldots,B_M(q_{49})\}$,
where $q_s$ is the left-lower corner vertex of $B_M(q_s)$, the same as before.
For example, we may think
$q_1=(0,0)$, $q_2=(1,0)$, $q_3=(0,1)$, $q_4=(-1,0)$, $q_5=(0,-1)\ldots.$
For each vertex $q_s$, we work on $\{B_{M}(q_s+(7i,7j))\}$ for all
integers $i$
and $j$.
In words, they are the $M$-square lattice on the plane at $7M$ apart.
With this definition,
\[
\bigcup_{s=1}^{49} \bigcup_{i,j} B_{M}
\bigl(q_s+(7i,7j)\bigr)=\mathbf{Z}^2.
\]
We also work on $\{G(B_M(q_s+(7i,7j)))\}$ for all $i$ and $j$.
By our definition, for $q_s$, these $7M$-squares
$\{G(B_M(q_s+(7i,7j)))\}$ are disjoint for all the different
$i$ or $j$ and the union of all these $7M$-squares is $\mathbf{Z}^2$.

For $n^{2/3} \leq m \leq n/2$ and $q_s$, we denote by $R_{M}(q_s, m,n)$
the number of squares of
$\{B_M(q_s+(7i,7j))\}$ that contain at least a vertex $v\in B(m)\cap
(S_M(z,\gamma_n)\cup D(z,\gamma_n))$
for all possible integers $i$ and $j$.
Note that for each $u\in B_M(q_s+(7i,7j))$, its $M$-bridge
%
%e4.8 ###
\begin{equation}\label{e4.7}
l_{u, v}\subset G\bigl(B_M\bigl(q_s+(7i,7j)\bigr)\bigr).
\end{equation}
Note also that
%
%e4.9 ###
\begin{equation}\label{e4.8}
\sum_{s=1}^{49} R_{M}(q_s,m,n)\geq
\bigl|B(m) \cap\bigl(S_M(z,\gamma_n)\cup D(z,\gamma_n)\bigr)\bigr|/M^2.
\end{equation}
If $m$ is not an integer, we may define
$ R_{M}(q_s,m,n)=R_M(q_s, \lfloor m \rfloor, n)$.
With Corollary \ref{cor1} and (\ref{e4.8}), we have the following corollary.
\begin{cor}\label{cor2}
Under the same hypotheses as Theorem \ref{thm1}, there exists
$C=C(F, z, M)$ such that
\[
E \Biggl( \sum_{s=1}^{49} R_{M}(q_s,m,n) \Biggr)\geq C m.
\]
\end{cor}

%s5 ###
\section[Proof of Theorem 2]{Proof of Theorem \protect\ref{thm2}}\label{s5}
\setcounter{equation}{-1}

Before the proof, we need to introduce a martingale inequality obtained by
Newman and Piza \cite{newPiz95}.
Let $U_1, U_2, \ldots$ be disjoint edge subsets of $\mathbf{Z}^2$. We will express
configuration $\omega$ for each $k$
as $(\omega_k, \hat{\omega}_k)$, where $\omega_k$ (resp. the edges in
$\hat{\omega}_k$) is the restriction of $\omega$ to $U_k$
(resp. the edges in $\mathbf{Z}^2\setminus U_k$). We also have, for each $k$,
disjoint events $D^+_k$ and $D_k^-$ in $\mathcal{F}(U_k)$,
where $\mathcal{F}(U_k)$, for each $k$, is the sigma-field generated by
$t(e)$ for $e\in U_k$. With these two events, let
%
%e5.0 ###
\begin{equation}\label{e5.0}
H_k(\omega)= {a_{0,n}^+}(\hat{\omega}_k)- {a_{0,n}^-}(\hat{\omega}_k),
\end{equation}
where
\[
{a_{0,n}^+}(\hat{\omega}_k)=\inf_{\omega_k \in D^+_k} a_{0,n}(\omega_k,
\hat{\omega}_k)
\quad \mbox{and}\quad
{ a_{0,n}^-}(\hat{\omega}_k)=\sup_{\omega_k
\in D^-_k} a_{0,n}(\omega_k, \hat{\omega}_k).
\]
Using these definitions, Newman and Piza \cite{newPiz95} proved in their
Theorem 8 the following lemma.
\begin{lem}[(Newman and Piza \cite{newPiz95})]\label{lem8}
If $U_k$, $D_k$, and $H_k$ satisfy the following:
\begin{longlist}
\item Conditional on $\mathcal{F}(\mathbf{Z}^2\setminus\bigcup_k U_k)$, then
$\mathcal{F}(U_i)$ and $\mathcal{F}(U_j)$
are mutually independent for $i \neq j$.
\item There exist positive $p$ and $q$ such that for any $k$
\[
P\bigl(\omega_k \in D_k^-|\mathcal{F}(\mathbf{Z}^2 \setminus U_k)\bigr) \geq p
\quad \mbox{and}\quad
P\bigl(\omega_k \in D_k^+ |\mathcal{F}(\mathbf{Z}^2 \setminus U_k)\bigr) \geq q
\qquad \mbox{a.s.}
\]
\item For every $k$, $H_k \geq0$ a.s.

Suppose that, for $\varepsilon>0$ and each $k$, $F_k\subset\mathcal{F}(\mathbf{Z}^2)$
is a subset of event $\{H_k\geq\varepsilon\}$.
Then
\[
\sigma^2( a_{0,n})\geq pq\varepsilon^2 \sum_{k} P(F_k)^2.
\]
\end{longlist}
\end{lem}

To apply Lemma \ref{lem8}, we set all vertices on $\mathbf{Z}^2$ in a spiral ordering
starting from the origin. We denote these vertices
by $\{(i_t, j_t)\}$ for $t=1,2,\ldots.$
Now we define vertex sets $U_1=G(B_M(q_1+(7i_1, 7j_1)))$,
$U_2=G(B_M(q_1+(7i_2,7j_2)))$, $U_3=G(B_M(q_1+(7i_3, 7j_3))),\ldots,$
which is a spiral ordering of these $7M$-squares.
Recall that our squares are the sets of vertices, but it is easy to reconsider
them as the edges in these squares
without the edges in the top and right boundaries.

Note that with this ordering, $U_1,\ldots, U_k, \ldots$ eventually
cover all $\mathbf{Z}^2$,
and
%
%e5.1 ###
\begin{equation}\label{e5.1}
B(m) \subset\bigcup_{k=1}^{m^2} U_k.
\end{equation}
Since $U_i\cap U_j=\varnothing$ for $i\neq j$, (i) in Lemma \ref{lem8} holds.
Let $D^-_k$ be the event that all edges in $U_k$ are $1$-edges and let $D^+_k$
be the event that
all edges in $U_k$ are $z^+$-edges. Since $U_k$ is finite, then
\[
P\bigl(\omega_k \in D_k^-|\mathcal{F}(\mathbf{Z}^2 \setminus U_k)\bigr)
\geq \vec{p}_c^{100M^2}
\]
and
\[
P\bigl(\omega_k \in D_k^+|\mathcal{F}(\mathbf{Z}^2 \setminus U_k)\bigr)
\geq \bigl[1-P\bigl(t(e)=1\bigr)- \delta_3\bigr]^{100M^2},
\]
where $\delta_3$ is defined in (\ref{e4.0}).

Therefore, (ii) in Lemma \ref{lem8} is satisfied if $\delta_3$ is small enough.
Note that $a_{0,n}$ is a coordinatewise nondecreasing function of
$\omega$, so (iii) holds.

Let $F_k(q_1) $ be the event that:
\begin{longlist}[(b)]
\item[(a)] $\gamma_n$, defined in Section \ref{s1},
has to use at least a $z^+$-edge of $U_k$ or
\item[(b)] there is an $M$-broken bridge $l_{u_i, v_i}\subset U_k$
($1\leq i\leq \tau$) for $\gamma_n$ such that $l_{u_i,v_i}$ contains at
least one $z^+$-edge.
\end{longlist}
We will show that
%
%e5.2 ###
\begin{equation}\label{e5.2}
F_k(q_1)\subset \bigl\{H_k\geq\min\{2, z-1\}\bigr\},
\end{equation}
so Lemma \ref{lem8} can be applied.

On (a), for $\omega=(\omega_k, \hat{\omega}_k)\in F_k(q_1)$, note
that if all $z^+$-edges in $U_k$
are changed to be 1-edges, the passage time $T(\gamma_n)$ is at least saved
by $z-1$, so
%
%e5.3 ###
\begin{equation}\label{e5.3}
\qquad
a_{0,n}(\omega) =T(\gamma_n)(\omega)\geq T(\gamma_n)(\omega^1_k,
\hat{\omega}_k)+(z-1) \geq a^-_{0,n}(\hat{\omega}_k)+(z-1),
\end{equation}
where $\omega^1_k$ is the configuration in $U_k$ such that all edges in $U_k$
have value one, and
$T(\gamma_n)(\omega^1_k, \hat{\omega}_k)$ is the passage time for path
$\gamma_n$, but with configuration $(\omega^1_k, \hat{\omega}_k)$.

On the other hand, we denote by $\omega^+_k$ the
configuration in $U_k$ such that
\[
a_{0,n}^+(\hat{\omega}_k)=\inf_{\omega_k \in D^+_k}
a_{0,n}(\omega_k,\hat{\omega}_k)=a_{0,n}(\omega_k^+, \hat{\omega}_k).
\]
For the configuration $\omega_k^+$, all edges in $U_k$ have values
larger than $z$.
With this new configuration $(\omega_k^+, \hat{\omega}_k)$, if an
optimal path for
$a_{0,n}(\omega_k^+, \hat{\omega}_k)$ never passes through~$U_k$, then
%
%e5.4 ###
\begin{equation}\label{e5.4}
a_{0,n}(\omega_k^+, \hat{\omega}_k)=a_{0,n}(\omega).
\end{equation}
By (\ref{e5.3}) and (\ref{e5.4}),
%
%e5.5 ###
\begin{equation}\label{e5.5}
a_{0,n}^+(\hat{\omega}_k) = a_{0,n}(\omega_k^+,
\hat{\omega}_k)=a_{0,n}(\omega)\geq a_{0,n}^-(\hat{\omega}_k)+z-1.
\end{equation}
Therefore, by (\ref{e5.5}) we have
%
%e5.6 ###
\begin{equation}\label{e5.6}
H_k\geq(z-1).
\end{equation}
If all optimal paths for $a_{0,n}(\omega_k^+, \hat{\omega}_k)$ have to pass
through $U_k$, we denote by $\gamma^+_n$
an optimal path for the configuration $(\omega_k^+, \hat{\omega}_k)$. Then
we reduce the value of the edges in $U_k\cap\gamma_n^+$ from $z$ to
$1$ to
have
%
%e5.7 ###
\begin{eqnarray}\label{e5.7}
a_{0,n}^+(\hat{\omega}_k) &=& a_{0,n}(\omega_k^+, \hat{\omega}_k)
\nonumber\\[-8pt]
\\[-8pt]
&\geq & T(\gamma^+_n)(\omega_k^1, \hat{\omega}_k)+(z-1)\geq a_{0,n}^-
(\hat{\omega}_k)+(z-1),
\nonumber
\end{eqnarray}
where
$T(\gamma_n^+)(\omega^1_k, \hat{\omega}_k)$ is the passage time for path
$\gamma_n^+$ with configuration $(\omega^1_k, \hat{\omega}_k)$.
Therefore, we still have (\ref{e5.6}).

On (b), for $\omega=(\omega_k, \hat{\omega}_k)\in F_k(q_1)$,
$l_{u_i,v_i}\subset U_k$ has to contain at least
one \mbox{$z^+$-edge}. If we change all edges in $U_k$ from $z>1$ to $1$, then all
the \mbox{$z^+$-edges} in $l_{u_i,v_i}$ are changed to be $1$-edges. If we go along
the bridge $l_{u_i,v_i}$ from $u_i$ to $v_i$,
we at least save time two, compared with going along $\gamma(u_i, v_i)$ from
$u_i$ to $v_i$. Therefore,
%
%e5.8 ###
\begin{equation}\label{e5.8}
a_{0,n}(\omega) \geq a^-_{0,n}(\hat{\omega}_k)+ 2.
\end{equation}
If an optimal path for
$a_{0,n}(\omega_k^+, \hat{\omega}_k)$ never passes through $U_k$, then
by (\ref{e5.8}),
%
%e5.9 ###
\begin{equation}\label{e5.9}
a_{0,n}(\omega_k^+, \hat{\omega}_k)=a_{0,n}(\omega)\geq a^-
_{0,n}(\hat{\omega}_k)+ 2.
\end{equation}
If all optimal paths for $a_{0,n}(\omega_k^+, \hat{\omega}_k)$ have to pass
through $U_k$,
by the same reason in (\ref{e5.7}), we have
%
%e5.10 ###
\begin{equation}\label{e5.10}
H_k\geq2.
\end{equation}
Together with (\ref{e5.9}) and (\ref{e5.10}), on (b), we also have
\[
H_k \geq\min(z-1, 2).
\]
Thus, (\ref{e5.2}) follows.
It follows from Lemma \ref{lem8}
that there exists $C=C(F, M, \delta_5, z)$ such that
%
%e5.11 ###
\begin{equation}\label{e5.11}
\sigma^2(a_{0,n})\geq C\sum_{k} [P(F_k(q_1))]^2.
\end{equation}

By Lemma 1 in Newman and Piza \cite{newPiz95}, we have
%
%e5.12 ###
\begin{equation}\label{e5.12}
\sigma^2(a_{0,n})\geq C (\log n )^{-1}
\Biggl(\sum_{m=1}^{n^2/4} m^{-3/2}
 \Biggl[ \sum_{k=1}^m P(F_k(q_1)) \Biggr] \Biggr)^2.
\end{equation}
By (\ref{e5.2}), we have for $n^{2/3}\leq m\leq n/2$,
%
%e5.13 ###
\begin{equation}\label{e5.13}
E R_M(q_1,m,n) \leq\sum_{k=1}^{m^2} P(F_k(q_1)).
\end{equation}
This shows that
%
%e5.14 ###
\begin{equation}\label{e5.14}
\sigma^2(a_{0,n})\geq C (\log n )^{-1}
\Biggl(\sum_{m=n^{4/3}}^{n^2/4} m^{-3/2}
\bigl[E R_M\bigl(q_1, \sqrt{m}, n\bigr) \bigr] \Biggr)^2.
\end{equation}

Similarly, we have the same inequalities corresponding to
$q_2,\ldots, q_{49}$, to have
%
%e5.15 ###
\begin{equation}\label{e5.15}
\sigma^2(a_{0,n})\geq C (\log n )^{-1}
\Biggl(\sum_{m=n^{4/3}}^{n^2/4} m^{-3/2}
\bigl[E R_M\bigl(q_t,\sqrt{m}, n\bigr) \bigr] \Biggr)^2.
\end{equation}
If we sum all $t$ from $t=1$ to $t=49$ together, by a standard
inequality ($2ab \leq a^2+b^2$ for positive $a$ and $b$)
we have
%
%e5.16 ###
\begin{equation}\label{e5.16}
\qquad \quad
49\sigma^2(a_{0,n}) % \\
\geq C (2^{49} \log n )^{-1}
\Biggl(\sum_{m=n^{4/3}}^{n^2/4} m^{-3/2}
\Biggl[ \sum_{t=1}^{49} E R_M\bigl(q_t, \sqrt{m}, n\bigr) \Biggr] \Biggr)^2.
\end{equation}

By using Corollary \ref{cor2}
for each $m\geq n^{2/3}$ in (\ref{e5.16}), we have
\[
\sigma^2(a_{0,n})\geq C \log n.
\]
So Theorem \ref{thm2} follows.

\section*{Acknowledgments}

The author acknowledges the referee's comments. In particular, he
gratefully for a referee's
detailed comments, a simple proof for Lemma~\ref{lem1}, and pointing out
inaccuracies in the proofs of
Lemma \ref{lem2}, the counting argument in (\ref{e3.14}), and unclear statement of
Lemma \ref{lem5},
which resulted in an improved exposition. The author also would like to
acknowledge M. Takei for his many valuable comments.

\printaddresses

\end{document}